\documentclass[11pt,a4paper, fullpage]{amsart}

\usepackage[centertags]{amsmath}
\usepackage{latexsym,amsthm,amssymb,verbatim}
\usepackage{epsf}
\usepackage[scanall]{psfrag}
\usepackage{graphicx, color}

\newtheorem{defin}{Definition}[section]
\newtheorem{prop-def}[defin]{Proposition-Definition}
\newtheorem{lem}[defin]{Lemma}
\newtheorem{thm}[defin]{Theorem}
\newtheorem{remark}[defin]{Remark}
\newtheorem{cor}[defin]{Corollary}

\newtheorem{ex}[defin]{Example}
\newtheorem{claim}{Claim}
\newcommand{\e}{ \hfill $\diamond$}

\begin{document}

\title[Stallings Foldings]{Stallings' Foldings and Subgroups of Amalgams of Finite Groups}%
\author{L.Markus-Epstein}
\footnote{Supported in part at the Technion by a fellowship of the
Israel Council for Higher Education. Partially supported by the
Focus Center, "Group Theoretic Methods in Algebraic Varieties",
Department of Mathematics, Bar-Ilan University, funded by the
Israel Academy of Science. Partially supported by the Emmy Noether
Research Center, Department of Mathematics, Bar-Ilan University.
Partially supported by the Doctoral Committee, Bar-Ilan
University.
}%
\address{Department of Mathematics \\
Technion \\
Haifa 32000, Israel}%
\email{epstin@math.biu.ac.il}%

\begin{abstract}

 In the 1980's Stallings \cite{stal} showed that every finitely
generated subgroup of a free group is canonically represented by a
finite minimal immersion of a bouquet of circles. In terms of the
theory of automata, this is a minimal finite inverse automaton.
This allows for the deep algorithmic theory of finite automata and
finite inverse monoids to be used to answer questions about
finitely generated subgroups of free groups.

In this paper we attempt to apply the same methods to other
classes of groups. A fundamental new problem is that the Stallings
folding algorithm must be modified to allow for ``sewing'' on
relations of non-free groups. We look at the class of groups that
are amalgams of finite groups. It is known that these groups are
locally quasiconvex and thus all finitely generated subgroups are
represented by finite automata. We present an algorithm to compute
such a finite automaton and use it to solve various algorithmic
problems.

\end{abstract}
\maketitle

\pagestyle{headings}
%

%%%%%%%%%%%%%%%%%%%%%%%%%%%%%%%%%%%%%%%%

\section{Introduction} \label{section:Introduction}

This paper has as its main goal the extension of results from the
case of studying subgroups of free groups to that of other classes
of finitely presented groups via  techniques in automata theory
and the theory of inverse semigroups. The main idea is that a
finitely generated subgroup $H$ of a ``nice'' group $G$ can be
represented by a finite directed graph labelled by generators of
$G$. From an automata theoretical point of view, the graph is a
finite inverse automaton, and from a topological point of view it
is an immersion over a bouquet of circles. This convergence of
ideas from group theory, topology, and the theory of finite
automata and finite semigroups allows for a rich interaction of
ideas and methods from many different fields.

In the free group, finitely generated subgroups correspond
precisely to finite inverse automata or equivalently to finite
immersions over a bouquet of circles. This object can be
constructed algorithmically by the process of Stallings foldings
\cite{stal}. It can be shown that every finitely generated
subgroup $H$ of a free group $FG(X)$ over  a set of generators $X$
corresponds to a uniquely determined such finite object
$\mathcal{A}(H)$ which is in fact a topological invariant of $H$.
Another important invariant of $H$ is its syntactic monoid
$\mathcal{M}(H)$ which is the transition monoid of
$\mathcal{A}(H)$. This is a finite inverse monoid. Thus
combinatorial and algorithmic properties of $H$ can be studied by
looking at the finite objects $\mathcal{A}(H)$ and
$\mathcal{M}(H)$. Since the theory of finite automata and finite
semigroups have rich algorithmic theories, non-trivial results can
be obtained in this way. In particular, this approach gives
polynomial time algorithms to solve the membership problem for $H$
(i.e. the Generalized Word Problem), the finite index problem, and
the computation of closures of $H$ in various profinite
topologies. On the other hand, the problem of checking purity,
that is, if $H$ is closed under taking roots, turns out to be
PSPACE-complete. See the articles \cite{a-o, b-m-m-w, kap-m,
mar_meak, m-s-w} for these and other examples of this approach.

%%%%%%%%%%%%%%%%%%%%%%%%%%%%%%%%%%%%%%%%%%%%%%%%%%%%%%%%%%%%%%%%%%%%%%

In general, the results mentioned above can not be extended to
every class of groups. That is because this immediately runs into
a problem: a theorem of Mikhailova \cite{l_s} shows that the
membership problem for the direct product of two non-Abelian free
groups is undecidable. Thus any hope of generalizing these results
to other classes of groups first must choose a class of groups
that are in some sense close to free groups, but far from the
direct product of two free groups!

%%%%%%%%%%%%%%%%%%%%%%%%%%%%%%%%%%%%%%%%%%%%%%%%%%%%%%%%%%%%%%%%%%%%%%

The groups considered in this paper   are \emph{amalgams of finite
groups}. As is well known, such groups are \emph{hyperbolic}
(\cite{b-f}) and \emph{locally quasiconvex} (\cite{ikap}). The
combination of these properties provides a fulfillment of the
above requirement.

%%%%%%%%%%%%%%%%%%%%%%%%%%%%%%%%%%%%%%%%%%%%%%%%%%%%%%%%%%%%%%%%%%%%%%%

Recall, that a group $G$ is \emph{locally quasiconvex} if and only
if every finitely generated subgroup $H$ of $G$ is quasiconvex. In
general, quasiconvexity of the subgroup $H$ depends on the
presentation of the group $G$. However  if the group $G$ is also
\emph{hyperbolic}, then the subgroup $H$ remains quasiconvex in
all finite presentations of $G$ (\cite{gro}).
%Thus for hyperbolic groups local quasiconvexity  is an intrinsic property of the
%group, rather than of a particular presentation.
This enables us to work with a fixed finite presentation of $G$
without loss of generality.

In \cite{gi_quas} Gitik  proved that the subgroup $H$ of the group
$G$ is quasiconvex if and only if the \emph{geodesic core}
$Core(G,H)$ of $Cayley(G,H)$ (which is the union of all closed
geodesics in the relative Cayley graph $Cayley(G,H)$ beginning at
the basepoint $H \cdot 1$) is finite. Thus local quasiconvexity of
the group $G$ (with a fixed finite presentation) ensures the
existence of a finite graph canonically associated with the given
subgroup $H$. Such a graph posses all the essential information
about the subgroup $H$ itself, therefore it can be used to study
properties of $H$.

%----------------------------------------------------------------

 However the geodesic core can not be constructed using a
generalization of Stallings' foldings algorithm. That is because
in amalgams, unlike in free groups, the theoretically well defined
notion of geodesic is ambiguous from computational and
constructible points of view. We are not familiar with any
rewriting procedure that computes geodesic words in  amalgams.
Since Stallings' foldings can be viewed as a simulation of a
rewriting procedure of freely reduced words in free groups, these
methods do not appear useful for a construction of geodesic cores.

In spite of this, normal (reduced) words do have a good
realization in amalgams given by their standard group
presentation. Indeed, there is a well known rewriting procedure
(\cite{l_s}) that given an element of an amalgam, computes its
normal (reduced) form. Such a rewriting is possible when the
amalgamated subgroup has a solvable membership problem in the
factors. Therefore it can be applied to elements of amalgams of
finite groups. This allows us to generalize Stallings' algorithm
following similar ideas and techniques.

Moreover, the following lemma of Gitik shows that geodesics and
strong normal paths are close to each other, which ensures in our
case the finiteness of the \emph{normal core} of $Cayley(G,H)$,
that is the union of all closed normal paths in $Cayley(G,H)$
starting at the basepoint $H \cdot 1$

\begin{lem}[Lemma 4.1 in \cite{gi_quas}] \label{lemma4.1}
If $G_1$ and $G_2$ are quasiconvex subgroups  of a hyperbolic
group  $G=G_1 \ast_{A} G_2$, then there exists a constant
$\epsilon \geq 0$ such that for any geodesic $\gamma \subset
Cayley(G)$ there exists a path $\gamma'$ in normal form with the
same endpoints as $\gamma$ with the following properties:
\begin{enumerate}
 \item [(1)] $\gamma \subset N_{\epsilon}(\gamma')$
\footnote{$ N_K(S)=\cup \{p \: | \: p \; {\rm is \; a \; geodesic
\; in} \; Cayley(G), \;  \iota(p) \in S, \; |p| \leq K \}$ is the
\emph{$K$-neighborhood} of $S$.}
  and   $\gamma' \subset N_{\epsilon}(\gamma)$,
 \item [(2)] an endpoint of any maximal monochromatic subpath of $\gamma'$ lies in $\gamma$ and is bichromatic in $\gamma$.
\end{enumerate}
\end{lem}

We explore normal cores and find that they can be defined  not
only theoretically, but constructively as well.
Theorem~\ref{ncore=reduced precover} says that the normal core of
$Cayley(G,H)$ is   a \emph{reduced precover} of $G$ (see
Definition \ref{def:reduced precover}), which is a restriction of
the notion of \emph{precovers}(see Section \ref{sec: precovers})
presented by Gitik in \cite{gi_sep}. Roughly speaking, one can
think of a reduced precover as a bunch of ``essential'' copies of
relative Cayley graphs of the free factors of $G$ glued to each
other according to the amalgamation. We prove (Corollary \ref{cor:
red precovers of the same subgroup are isomorphic}) that reduced
precovers determining the same subgroup are isomorphic.
Furthermore, our Maim Theorem (Theorem \ref{thm: unique reduced
precover}) states that given a finitely generated subgroup $H$ of
an amalgam  $G=G_1 \ast_A G_2$ there exists a unique reduced
precover determining $H$, which is the normal core of
$Cayley(G,H)$.

This constructive characterization of normal cores enables us to
present a quadratic algorithm (see Section \ref{subsec:
TheAlgorithm}) that given a finite set of subgroup generators of
$H$ constructs the normal core of $Cayley(G,H)$, where $G$ is an
amalgam of finite groups. Theorem~\ref{thm amalgam of finite grps}
provides the validity and the finiteness of the construction.

%%%%%%%%%%%%%%%%%%%%%%%%%%%%%%%%%%%%%%%%%%%%%%%%%%%%%%%%%%%%%%%%%%%

Thus the normal core $\Delta$ of $Cayley(G,H)$ posses properties
analogous to those of graphs constructed by the Stallings'
algorithm for finitely generated subgroups of free groups.
Geometrically, it can be viewed as the 1-skeleton of a topological
core of the covering space corresponding to $H$ of the standard
2-complex of $G$. Algebraically, $\Delta$ is an ``essential part''
of the relative Cayley graph $Cayley(G,H)$, and from the automata
theoretic point of view, it is a minimal finite inverse automaton
$\mathcal{A}$
%over the alphabet $X_1^{\pm} \cup X_2^{\pm}$
such that $L(\mathcal{A})=_G H$.

%%%%%%%%%%%%%%%%%%%%%%%%%%%%%%%%%%%%%%%%%%%%%%%%%%%%%%%%%%%%%%%%%%%%

Furthermore, Theorem~\ref{thm amalgam of finite grps} ensures the
canonicity of our construction, that is its independence from the
choice of subgroup generators, and guarantees that the resulting
graph $\Gamma(H)$ ``accepts'' all normal forms of elements from
$H$. We get the following corollary which gives an immediate
solution for the membership problem of $H$.
\begin{cor} \label{remark: normal elements in H}
A  normal word $g$ is in $H$ if and only if it labels a closed
path from $v_0$ to itself in $\Gamma(H)$.
\end{cor}
An application of normal cores  yields polynomial (mostly
quadratic) solutions for a nice list of algorithmic problems
concerning finitely generated subgroups of amalgams of finite
groups: the membership problem, the finite index problem, the
freeness problem, the power problem, the conjugacy problem, the
normality and the malnormality problems. Furthermore, the
separability problem can be solved in some particular cases and an
effective Kurosh decomposition for finitely generated subgroups in
the case of free products can be found. All these results are
presented in the PhD thesis of the author \cite{m-thesis}.

The present paper includes only the solution for the membership
problem as a demonstration of the effectiveness of our methods.
The rest of the above algorithmic problems and their solutions
will appear in our future papers \cite{m-algorithms, m-kurosh}.

%%%%%%%%%%%%%%%%%%%%%%%%%%%%%%%%%%%%%%%%%%%%%%%%%%%%%%%%%%%%%%%%%%%

Finally, we notice  that there are several generalization results
of Stallings' algorithm to other classes of groups. Schupp in
\cite{schupp} presents an algorithm for certain Coxeter groups and
surface groups of an extra-large type. Kapovich and Schupp
\cite{kap-s} make use of modified Stallings' foldings and the
minimization technique of Arzhantseva and Ol'shanskii \cite{a-o}
to present finitely generated subgroups of Coxeter groups and
Artin groups of extra-large type  and also of one-relator groups
with torsion by labelled graphs. Kapovich, Weidman, and Miasnikov
in \cite{kap-w-m} develop a combinatorial treatment of Stallings'
foldings  in the context of graphs of groups through the use of
the Bass-Serre theory. McCammond and Wise \cite{m_w} generalize
Stallings' algorithm for the class of coherence groups, however
the resulting graphs  are not canonical (they depend on the choice
of subgroup generators). Hence they are not suitable for solving
algorithmic problems for subgroups via their graphs.
Recently Miasnikov, Remeslennikov and Serbin have generalized
Stallings' algorithm to the class of fully residually free groups
\cite{m-r-s}. The developed methods were  applied  to solve a
collection of algorithmic problems concerning this class of groups
in \cite{k-m-r-s}.

\subsection*{Other Methods}  \
There have been a number of papers, where methods, not based on
Stallings' foldings, have been presented. One can use these
methods to treat finitely generated subgroups of amalgams of
finite groups. A topological approach can be found in works of
Bogopolskii \cite{b1, b2}. For the automata theoretic approach,
see papers of Holt and Hurt \cite{holt-decision, holt-hurt},
papers of Cremanns, Kuhn, Madlener and Otto \cite{c-otto,
k-m-otto}, as well as the recent paper of Lohrey and Senizergues
\cite{l-s}.

However the methods for treating finitely generated subgroups
presented in the above papers were applied to some particular
subgroup property. No one of these papers have as its goal a
solution of various algorithmic problems, which we consider as our
primary aim. We view the current paper as the first step in its
achieving. Similarly to the case of free groups (see
\cite{kap-m}), our combinatorial approach seems to be the most
natural one for this purpose. It yields reach algorithmic results,
as appear in our future papers \cite{m-algorithms, m-kurosh}.

%%%%%%%%%%%%%%%%%%%%%%%%%%%%%%%%%%%

\section{Acknowledgments}

I wish to deeply thank to my PhD advisor Prof. Stuart W. Margolis
for introducing me to this subject, for his help  and
encouragement throughout  my work on the thesis. I owe gratitude
to Prof. Arye Juhasz for his suggestions and many useful comments
during the writing of this paper. I gratefully acknowledge a
partial support at the Technion by a fellowship of the Israel
Council for Higher Education.

%%%%%%%%%%%%%%%%%%%%%%%%%%%%%%%%%%%%

\section{Labelled Graphs} \label{sec:LabeledGraphs}

The primary idea of this paper is to study finitely generated
subgroups of amalgams of finite groups by constructing subgroup
graphs exactly as in the case of free group. Hence we begin by
fixing the notation on graphs that will be used along this work.
In doing so we follow the notations used by Stallings in
\cite{stal} and Gitik in \cite{gi_sep}.

At the end of the section  we recall the notion of Stallings'
foldings and introduce a new graph operation which is an immediate
generalization of foldings for a non free group $G$. We prove that
both operations when applied to a subgroup graph $\Gamma(H)$, $H
\leq G$, do not affect the subgroup $H$.

%%%%%%%%%%%%%%%%%%%%%%%%%%%%%%%%%%%%%%%%%%%%%%%%%%%%%%%%%%%%%%%%%%%%%%%%%%%%%%%%%%%%%%%%%%%%%

%The following notations are taken from \cite{gi_sep} and \cite{gi_doub}.
%
A graph $\Gamma$ consists of two sets $E(\Gamma)$ and $V(\Gamma)$,
and two functions $E(\Gamma)\rightarrow E(\Gamma)$  and
$E(\Gamma)\rightarrow V(\Gamma)$: for each $e \in E$ there is an
element $\overline{e} \in E(\Gamma)$ and an element $\iota(e) \in
V(\Gamma)$, such that $\overline{\overline{e}}=e$ and
$\overline{e} \neq e$. The elements of $E(\Gamma)$ are called
\textit{edges}, and an $e \in E(\Gamma)$ is a \emph{directed edge}
of $\Gamma$, $\overline{e}$ is the \emph{reverse (inverse) edge}
of $e$. The elements of $V(\Gamma)$ are called \textit{vertices},
$\iota(e)$ is the \emph{initial vertex} of $e$, and
$\tau(e)=\iota(\overline{e})$ is the \emph{terminal vertex} of
$e$. We call them the \emph{endpoints} of the edge $e$.

\begin{remark}
{\rm A \emph{subgraph} of $\Gamma$ is a graph $C$  such that $V(C)
\subseteq V(\Gamma)$ and $E(C) \subseteq E(\Gamma)$. In this case,
by abuse of language, we write $C\subseteq \Gamma$.

Similarly, whenever we write $\Gamma_1 \cup \Gamma_2$ or $\Gamma_1
\cap \Gamma_2$  we always mean that the set operations are, in
fact,  applied to the vertex sets and the edge sets of the
corresponding graphs.

\e}
\end{remark}

A \emph{labelling} of $\Gamma$ by the set $X^{\pm}$ is a function
$$lab: \: E(\Gamma)\rightarrow X^{\pm}$$ such that for each $e \in
E(\Gamma)$, $lab(\overline{e})=(lab(e))^{-1}$.

The last equality enables one, when representing the labelled
graph $\Gamma$ as a directed diagram,  to represent only
$X$-labelled edges, because $X^{-1}$-labelled edges can be deduced
immediately from them.

A graph with a labelling function is called a \emph{labelled (with
$X^{\pm}$) graph}. A labelled graph is called \emph{well-labelled}
if
$$\iota(e_1)=\iota(e_2), \; lab(e_1)=lab(e_2)\ \Rightarrow \
e_1=e_2,$$ for each pair of edges $e_1, e_2 \in E(\Gamma)$. See
Figure \ref{fig: labelled, well-labelled graphs}.

\begin{figure}[!h]
\psfrag{a }[][]{$a$} \psfrag{b }[][]{$b$} \psfrag{c }[][]{$c$}
\psfrag{e }[][]{$e_1$}
\psfrag{f }[][]{$e_2$}
\psfragscanon \psfrag{G }[][]{{\Large $\Gamma_1$}}
\psfragscanon \psfrag{H }[][]{{\Large $\Gamma_2$}}
\psfragscanon \psfrag{K }[][]{{\Large $\Gamma_3$}}
\includegraphics[width=\textwidth]{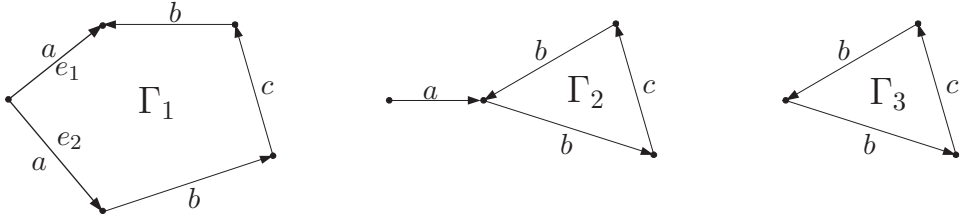}
\caption[The construction of $\Gamma(H_1)$]{ \footnotesize {The
graph $\Gamma_1$ is labelled with $\{a,b,c\} ^{\pm}$, but it is
not well-labelled. The graphs $\Gamma_2$ and $\Gamma_3$ are
well-labelled with $\{a,b,c\} ^{\pm}$.}
 \label{fig: labelled, well-labelled graphs}}
\end{figure}

 The label of a path $p=e_1e_2 \cdots
e_n$ in $\Gamma$, where $e_i \in E(\Gamma)$, is the word $$lab(p)
\equiv lab(e_1)\cdots lab(e_n) \in (X^{\pm})^*.$$ Notice that the
label of the empty path is the empty word. As usual, we identify
the word $lab(p)$ with the corresponding element in $G=gp\langle X
\; | \; R\rangle$.

Note that elements of $G=gp \langle X |R \rangle$ are equivalence
classes of words. However it is customary to blur the distinction
between a word $u$ and the equivalence class containing $u$. We
will distinguish between  them by using different equality signs:
\fbox{``$\equiv$''}
\footnote{Throughout the present paper boxes are used to emphasize
the notation.}
for the equality of two words and \fbox{``$=_G$''} to denote the
equality of two elements of $G$, that is the equality of two
equivalence classes.

A path $p=e_1e_2 \cdots e_n$ is freely reduced if $e_{i+1} \neq
\overline{e_i}$ for all $1 \leq i \leq n-1$.

\begin{remark}
{\rm If $\Gamma$ is a well-labelled graph then a path $p$ in
$\Gamma$ is freely reduced if and only if $lab(p)$ is a freely
reduced word.} \e
\end{remark}

Denote the pair consisting  of the graph $\Gamma$ and the
\emph{basepoint} (a distinguished vertex of the graph $\Gamma$)
$v_0$ by $(\Gamma, v_0)$ and call it a \emph{pointed graph}.

Following the notation of Gitik, \cite{gi_sep}, we denote the set
of all closed paths in $\Gamma$ starting at $v_0$  by
\fbox{$Loop(\Gamma, v_0)$},  and the image of $lab(Loop(\Gamma,
v_0))$ in $G$  by \fbox{$Lab(\Gamma, v_0)$}. More precisely,
$$Loop(\Gamma, v_0)=\{ p \; | \; p {\rm  \ is \ a \ path \ in \ \Gamma \
with} \ \iota(p)=\tau(p)=v_0\}, $$
$$Lab(\Gamma,v_0)=\{g \in G \; | \; \exists p \in Loop(\Gamma,
v_0) \; : \; lab(p)=_G g \}.$$

\begin{remark}[\cite{gi_sep}]
{\rm It is easy to see that $Lab(\Gamma, v_0)$ is a subgroup of
$G$. } \e
\end{remark}

\begin{remark} \label{remark: trivial subgroup}
{\rm If $V(\Gamma)=\{v_0\}$ and $E(\Gamma)=\emptyset$ then we
assume that $Lab(\Gamma, v_0)=\{1\}$. } \e
\end{remark}

\begin{remark}
{\rm We say that $H=Lab(\Gamma, v_0)$ is the subgroup of $G$
\emph{determined} by the graph $\Gamma$. Thus any pointed graph
labelled by $X^{\pm}$, where $X$ is a generating set of the group
$G$, determines a subgroup of $G$. This argues the use of the name
\emph{subgroup graphs} for such graphs. } \e
\end{remark}

%-----------------------------------------------------------------

%-------------Inverse Automata and Immersions---------------------

As is well known \cite{b-m-m-w, mar_meak, kap-m}, well-labelled
graphs, presented above combinatorially, can be viewed as
algebraical, topological, geometrical and automata-theoretical
objects as well.   The detailed exploration of various connections
between combinatorial group theory, semigroup theory and formal
language theory can be found in \cite{d-margolis-s}.

Thus a finite pointed graph $(\Gamma,v_0)$ well-labelled with
$X^{\pm}$ can be viewed as the inverse automaton  with the same
initial-terminal state $v_0$:
 $$\mathcal{A}=(V(\Gamma),X^{\pm},\delta,v_0,\{v_0\}),$$
where $\delta:V(\Gamma) \times X^{\pm} \rightarrow V(\Gamma)$,
usually denoted $\delta(v,x)=v \cdot x$, satisfies $\delta(v,x)=w$
if and only if there exist $e \in E(\Gamma)$ such that
$\iota(e)=v$, $\tau(e)=w$ and $lab(e) \equiv x$.   The
representation of $(\Gamma,v_0)$ is the positive state graph of
$\mathcal{A}$ and $L(\mathcal{A})=lab(Loop(\Gamma, v_0))$. The
reader is referred to \cite{b-m-m-w} for the missing definitions.

As usual, $\delta$ is extended to a (partial) function on
$V(\Gamma) \times (X^{\pm})^*$ by letting $v \cdot 1=v$ and $v
\cdot(ua)=(v \cdot u) \cdot a$ (if this is defined) for all $v \in
V(\Gamma)$, $u \in (X^{\pm})^*$ and $a \in X^{\pm}$.  Thus if $v,
w \in V(\Gamma)$ and $p$ is a path in $\Gamma$ such that
$$\iota(p)=v, \ \tau(p)=w \ {\rm and } \ lab(p)\equiv u,$$
then, following the automata theoretic notation, we simply write
$v \cdot u=w$ to summarize this situation.

By abuse of language, we say that a word $w$ is \emph{accepted} by
the graph $(\Gamma,v_0)$ if and only if there exists a path $p$ in
$\Gamma$ closed at $v_0$, $\iota(p)=\tau(p)=v_0$ such that $lab(p)
\equiv w$, that is $v_0 \cdot w=v_0$.

%%%%%%%%%%%%%%%%%%%%%%%%%%%%%%%%%%%%%%%%%%%%%%%%%%%%%%%%%%%%%%%%%%%%%%

\subsection*{Morphisms of Labelled Graphs}
\label{sec:Morphisms Of Well-Labelled Graphs}

%\begin{defin} \label{def: morphism lab graphs}
Let $\Gamma$ and $\Delta$ be graphs labelled with $X^{\pm}$. The
map $\pi:\Gamma \rightarrow \Delta$ is called a \emph{morphism of
labelled graphs}, if $\pi$ takes vertices to vertices, edges to
edges, preserves labels of directed edges and has the property
that
$$ \iota(\pi(e))=\pi(\iota(e)) \ {\rm and } \
\tau(\pi(e))=\pi(\tau(e)), \ \forall e\in E(\Gamma).$$
An injective morphism of labelled graphs is called an
\emph{embedding}. If $\pi$ is an embedding then we say that the
graph $\Gamma$ \emph{embeds} in the graph $\Delta$.
%\end{defin}

%-----------------------------------------------------------------------------

A \emph{morphism of pointed labelled graphs} $\pi:(\Gamma_1,v_1)
\rightarrow (\Gamma_2,v_2)$  is a morphism of underlying labelled
graphs $ \pi: \Gamma_1\rightarrow \Gamma_2$ which preserves the
basepoint $\pi(v_1)=v_2$. If $\Gamma_2$ is well-labelled then
there exists at most one such morphism (\cite{kap-m}).

%-----------------------------------------------------------------

\begin{remark}[\cite{kap-m}] \label{unique isomorphism}
{\rm  If two pointed well-labelled (with $X^{\pm}$) graphs
$(\Gamma_1,v_1)$ and $(\Gamma_2,v_2)$  are isomorphic, then there
exists a unique isomorphism $\pi:(\Gamma_1,v_1) \rightarrow
(\Gamma_2,v_2)$. Therefore $(\Gamma_1,v_1)$ and $(\Gamma_2,v_2)$
can be identified via $\pi$. In this case we sometimes write
$(\Gamma_1,v_1)=(\Gamma_2,v_2)$.} \e
\end{remark}

%%%%%%%%%%%%%%%%%%%%%%%%%%%%%%%%%%%%%%%%%%%%%%%%%%%%%%%%%%%%%%%%%%%%%%

\subsection*{Graph Operations} \label{sec:SubgroupGraphs}

Recall that a \emph{Stallings' folding} is an identification of a
pair of distinct edges with the same initial vertex and the same
label. The operation of  \emph{``cutting hairs''} consists of
removing from the graph edges whose terminal vertex has  degree
$1$ (see Figure \ref{fig: labelled, well-labelled graphs}: the
graph $\Gamma_2$ is obtained from the graph $\Gamma_1$ by folding
the edges $e_1$ and $e_2$; the graph $\Gamma_3$ is obtained from
the graph $\Gamma_2$ by cutting the hair edge labelled by $a$). As
is well known \cite{stal, mar_meak, kap-m}, these graph operations
don't affect the corresponding subgroup of a free group. The
following lemma demonstrates the similar behavior in the case of
finitely presented non free groups.

\begin{lem} \label{foldings-cutting-hairs}
    Let $G=gp\langle X|R \rangle$ be a finitely presented group. Let
    $\Gamma$ be a graph labelled with $X^{\pm}$ and let $\Gamma'$ be a
    graph labelled with $X^{\pm}$ obtained from $\Gamma$ by a single
    folding or by ``cutting'' a single hair. Then $Lab(\Gamma,v_0)=Lab(\Gamma',v_0')$, where
    $v_0$ is the basepoint of $\Gamma$ and $v_0'$ is the corresponding
    basepoint of $\Gamma'$.
\end{lem}
\begin{proof}
%\textbf{Proof.}
%
Let $F(X)$ be a free group with  finite free basis $X$. Let
$$lab: \: E(\Gamma)\rightarrow X^{\pm}$$
be the labelling function of $\Gamma$. The function $lab$ extends
to the labelling of paths of $\Gamma$ such that the label of a
path $p=e_1e_2 \cdots e_n$ in $\Gamma$, is the word $lab(p) \equiv
lab(e_1)\cdots lab(e_n) \in (X^{\pm})^*$. Denote by
$Lab_{F(X)}(\Gamma)$ the image of $lab(Loop(\Gamma))$ in $F(X)$.

 As is well known, \cite{stal, mar_meak, kap-m},  foldings and cutting hairs don't
affect the fundamental group of the graph, i.e.
$$Lab_{F(X)}(\Gamma,v_0)=Lab_{F(X)}(\Gamma',v_0').$$
Since the homomorphism $(X^{\pm})^* \rightarrow G$ factors through
$F(X)$
 $$(X^{\pm})^* \rightarrow F(X) \rightarrow G,$$
we conclude that $Lab(\Gamma,v_0)=Lab(\Gamma',v_0')$.

\end{proof}

%---------------------------------------------------------------------------
Let $f_1$ and $f_2$ be a pair of folded edges of the graph
$\Gamma$ with labels $x$ and $x^{-1}$, respectively. Hence the
path $\overline{f_1}f_2$ in $\Gamma$ is labelled by the trivial
relator $x^{-1}x$. The folding operation applied to the edges
$f_1$ and $f_2$ implies the identification of the endpoints of
$\overline{f_1}f_2$. Thus the natural extension of such operation
to the case of a non free group $G$ is an identification of the
endpoints of paths labelled by a relator.

\begin{figure}[!h]
\psfrag{a }[][]{$a$} \psfrag{b }[][]{$b$} \psfrag{c }[][]{$c$}
\psfrag{d }[][]{$d$} \psfrag{p }[][]{$p$}
\psfrag{v }[][]{$v_1$}
\psfrag{u }[][]{$v_2$}
\psfrag{w }[][]{$\vartheta$}
\psfragscanon \psfrag{G }[][]{{\Large $\Gamma_1$}}
\psfragscanon \psfrag{H }[][]{{\Large $\Gamma_2$}}
\includegraphics[width=\textwidth]{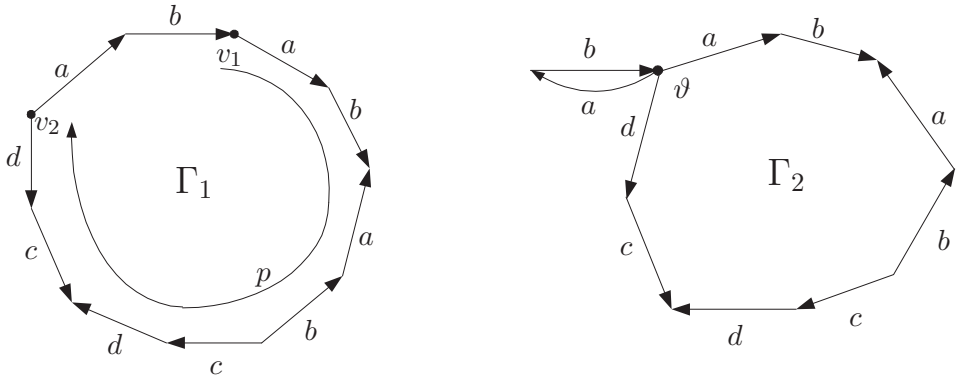}
\caption{ \footnotesize {Let $G=gp\langle a,b,c,d \; | \;
aba^{-1}b^{-1}cdc^{-1}d^{-1} \rangle$. Thus the graph $\Gamma_2$
is obtained from the graphs $\Gamma_1$ by the identification of
the endpoints $v_1$ and $v_2$ of the path $p$ labelled by the
defining relator of $G$. Note that the resulting graph $\Gamma_2$
is not well-labelled (at $\vartheta$).}
 \label{fig:IdentificationRelator}}
\end{figure}

\begin{defin} \label{def: identificacion-def-relator}
{\rm Let $\Gamma$ be a  graph labelled with $X^{\pm}$. Suppose
that $p$ is a path of $\Gamma$  with
$$v_1= \iota(p)\neq \tau(p)=v_2 \ {\rm and } \ lab(p) =_G 1.$$
Let $\Delta$ be a graph labelled with $X^{\pm}$ defined as
follows.

The vertex set of $\Delta$ is a vertex set of $\Gamma$ with
$\iota(p)$ and $\tau(p)$ removed and a new vertex $\vartheta$
added (we think of the vertices $\iota(p)$ and $\tau(p)$ as being
identified to produce vertex $\vartheta$):
$$V(\Delta)=(V(\Gamma) \setminus \{\iota(p),\tau(p)\}) \cup
\{\vartheta\}.$$
The edge set of $\Delta$ is the edge set of $\Gamma$:
$$E(\Delta)=E(\Gamma).$$ The endpoints and arrows for the edges of $\Delta$ are
defined in a natural way. Namely, if $e \in E(\Delta)$ and
$\iota(e), \tau(e) \not\in \{v_1,v_2\}$  then we put
$\iota_{\Delta}(e)=\iota_{\Gamma}(e)$. Otherwise
$\iota_{\Delta}(e)=\vartheta$ if
  $\iota_{\Gamma}(e) \in \{v_1,v_2\}$ and
  $\tau_{\Delta}(e)=\vartheta$ if $\tau_{\Gamma}(e)\in
  \{v_1,v_2\}$.

We define labels on the edges of $\Delta$ as follows:
$lab_{\Delta}(e) \equiv lab_{\Gamma}(e)$ for all $e \in
E(\Gamma)=E(\Delta)$.

Thus $\Delta$ is a graph labelled with $X^{\pm}$. In this
situation we say that $\Delta$ is obtained from $\Gamma$ by the
\emph{identification of a relator.} See Figure
\ref{fig:IdentificationRelator}}

\end{defin}

\begin{remark} \label{remark: identif-relator}
{\rm There exists an  epimorphism of pointed labelled graphs
$\phi: (\Gamma, v_0) \rightarrow (\Delta, u_0)$
 such that
$$\phi(v)=\left\{%
\begin{array}{ll}
   v, & \hbox{$v \not\in \{v_1,v_2\};$} \\
    \vartheta, & \hbox{otherwise.} \\
\end{array}%
\right.$$

Thus $u_0=\phi(v_0)$ and paths  in $\Delta$ are images of paths in
$\Gamma$. However, in order to simplify the notation we omit the
use of the graph morphism $\phi$. We  say that $\alpha'$ is a
vertex/edge/path in $\Delta$  \emph{corresponding } to the
vertex/edge/path $\alpha$ in $\Gamma$, instead of saying that
$\alpha'=\phi(\alpha)$ is the \emph{image} of $\alpha$ in
$\Delta$. We treat $\Delta$ as a graph constructed from $\Gamma$
in the combinatorial way described in Definition \ref{def:
identificacion-def-relator}.

}\e
\end{remark}

%---------------------------------------------------------------------------
\begin{lem} \label{identificacion-def-relator}
    Let $G=gp\langle X|R \rangle$ be a finitely presented group. Let
    $\Gamma$ be a graph well-labelled with $X^{\pm}$. Let $p$ be a freely reduced path in $\Gamma$
     with $lab(p)=_G 1$ such that  $\iota(p) \neq \tau(p)$.

    Let $\Gamma'$ be a
    graph obtained from $\Gamma$ by the
    identification of the endpoints of $p$. Then $Lab(\Gamma,v_0)=Lab(\Gamma',v_0')$, where
    $v_0$ is the basepoint of $\Gamma$ and $v_0'$ is the corresponding
    basepoint
    %
%%%%%
%\footnote{That is $v_0'=v_0$, if $v_0 \not\in \{\iota(p), \tau(p)
%\}$, and $v_0'=\vartheta$ (see Definition \ref{def:
%identificacion-def-relator}), otherwise.}
%%%%%
    %
     of $\Gamma'$.
\end{lem}
%
%\textbf{Proof }:
\begin{proof}
Let  $q \in Loop(\Gamma,v_0)$. The identification of the endpoints
of the path $p$ keeps closed paths of $\Gamma$ closed (because the
graph morphism $\phi: (\Gamma,v_0) \rightarrow (\Gamma',v_0')$,
see Remark \ref{remark: identif-relator}, preserves endpoints).
Thus the path $q'$ in $\Gamma'$ corresponding to the path $q$ in
$\Gamma$ (that is obtained from $q$ by the identification of the
endpoints of $p$) is closed at $v_0'$ if $p$ is a subpath of $q$
or if it is not a subpath of $q$. Thus $Loop(\Gamma,v_0) \subseteq
Loop(\Gamma',v'_0)$. Hence $ Lab(\Gamma, v_0) \subseteq
Lab(\Gamma',v'_0)$.

%Let $u \in Lab(\Gamma,v_0)$. Then there is $q \in
%Loop(\Gamma,v_0)$ such that $lab(q)=_G u$.
%
%Then either $p$ is or is not a subpath of $q$. In the second case
%the identification of the endpoints of $p$ doesn't affect the path
%$q$. Therefore $q \in Loop(\Gamma',v'_0)$ and hence $lab(q)=_G u
%\in Lab(\Gamma',v'_0)$. If $p$ is a subpath of $q$ then
%$q=q_1pq_2$. The desired identification remains the corresponding
%(to $q$) path $q'$ in $\Gamma'$ closed. Thus $q' \in
%Loop(\Gamma',v'_0)$. Moreover
%$lab(q')=lab(q_1)lab(p)lab(q_2)=lab(q)=_G u \in
%Lab(\Gamma',v'_0)$. Thus $ Lab(\Gamma, v_0) \subseteq
%Lab(\Gamma',v'_0)$.

Suppose now that  $w \in Lab(\Gamma',v'_0)$. Then there is $q' \in
Loop(\Gamma',v'_0)$ such that $lab(q')=_G w$. If $q'$ exists in
$\Gamma$ (i.e. $q' \in Loop(\Gamma,v_0) \cap Loop(\Gamma',v'_0)$)
then $w=_G lab(q') \in Lab(\Gamma,v_0)$.

Otherwise, $q' \in Loop(\Gamma',v'_0) \setminus Loop(\Gamma,v_0)$.
Let $p'$ be the path corresponding to the path $p$ in $\Gamma'$
and $\vartheta \in V(\Gamma')$ be the vertex corresponding to the
identified endpoints of the path $p$. Thus
$$\vartheta=\iota(p')=\tau(p'), \ lab(p) \equiv lab(p').$$
Hence  the following is possible.
\begin{itemize}
  \item
    $p'$ is not a subpath of $q'$.

    Then there is a decomposition $q'=q'_1q'_2 \ldots q'_k$ such
    that
    $$\iota(q'_1)=\tau(q'_k)=v'_0, \ \tau(q'_i)=\iota(q'_{i+1})=\vartheta, 1 \leq i \leq k-1,$$
    where
    $q'_i$ is a path in $\Gamma \cap \Gamma'$
    and $q'_iq'_{i+1}$ is a path in $\Gamma'$ which doesn't exist in $\Gamma$ (see
Figure \ref{fig:ProofIdentRelator}).
    It means that $\tau(q'_i)$ and $\iota(q'_{i+1})$ are different endpoints of the path $p$ in
$\Gamma$.
\begin{figure}[!h]
\psfrag{a }[][]{$q'_i$} \psfrag{b }[][]{$q'_{i+1}$} \psfrag{p
}[][]{$p$} \psfrag{q }[][]{$p'$}
\psfrag{v }[][]{$v_1$}
\psfrag{u }[][]{$v_2$}
\psfrag{w }[][]{$\vartheta$}
\psfragscanon \psfrag{G }[][]{{\Large $\Gamma$}}
\psfragscanon \psfrag{H }[][]{{\Large $\Gamma'$}}
\includegraphics[width=\textwidth]{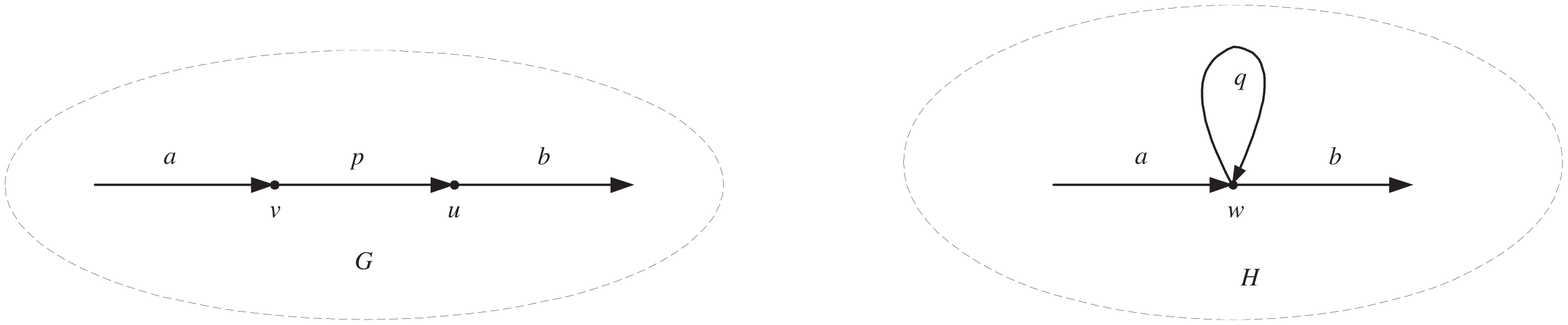}
\caption{ \label{fig:ProofIdentRelator}}
\end{figure}

Hence if $k=1$ then the path $q'$ is in $\Gamma \cap \Gamma'$.
Therefore $$w=_G lab(q') \in Lab(\Gamma, v_0).$$

Otherwise, let
$$p_i=\left\{%
\begin{array}{ll}
    p, & \hbox{$\tau(q'_i)=\iota(p)$;} \\
    \overline{p}, & \hbox{$\tau(q'_i)=\tau(p)$.} \\
\end{array}%
\right.$$
Thus $ q=q'_1p_1q'_2p_2 \ldots p_{k-1}q'_k$ is a path in $\Gamma$
closed at $v_0$. Since $lab(p_i)=_G 1$, we have
\begin{eqnarray}
 lab(q) & \equiv & lab(q'_1)lab(p_1)lab(q'_2)lab(p_2) \ldots lab(p_{k-1})lab(q'_k) \nonumber \\
  & =_G & lab(q'_1)lab(q'_2) \ldots lab(q'_k) \equiv lab(q'). \nonumber
\end{eqnarray}
Thus $w=_G lab(q')=_G lab(q) \in Lab(\Gamma, v_0)$. \\

  \item
    $p'$ is a subpath of $q'$.

    The occurrences of $p'$ subdivide $q'$ into a concatenation of
    paths of the form
    $q'=q'_1p'_1q'_2p'_2 \ldots p'_{k-1} q'_k$, where $ p'_i \in \{p', \overline{p'} \}$
    and the paths $q'_i$ don't involve $p$.

For each $1 \leq i \leq k$, the path $q'_i$   can be written as a
decomposition of subpaths in $\Gamma \cap \Gamma'$, and the
technique presented above (in the previous case) can be applied to
it. Hence for all $1 \leq i \leq k$, there exists a path $q_i
\subseteq \Gamma$ such that $\iota(q_i)=\iota(q'_i)$,
$\tau(q_i)=\tau(q'_i)$ and  $lab(q_i)=_G lab(q'_i)$.

Let $$p_i=\left\{%
\begin{array}{ll}
    p, & \hbox{$p'_i=p'$;} \\
    \overline{p}, & \hbox{$p'_i=\overline{p'}$.} \\
\end{array}%
\right.$$ Then $ q=q_1p_1q_2p_2 \ldots p_{k-1}q_k$ is  a path in
$\Gamma$ closed at $v_0$. Moreover,
\begin{eqnarray}
 lab(q) & \equiv & lab(q_1)lab(p_1)lab(q_2)lab(p_2) \ldots lab(p_{k-1})lab(q_k) \nonumber \\
  & =_G & lab(q'_1)lab(p'_1)lab(q'_2)lab(p'_2) \ldots lab(p'_{k-1})lab(q'_k)  \equiv lab(q'). \nonumber
\end{eqnarray}
\end{itemize}

Therefore $Lab(\Gamma,v_0)=Lab(\Gamma',v_0')$.

\end{proof}

%%%%%%%%%%%%%%%%%%%%%%%%%%%%%%%%%%%%%%%%%%%%%%%%%%%%%%%%%%%%%%%%%%%%%%%%%%%%%%%%%%%%%%%%%%%%%%%

\section{Subgroups and Covers} \label{sec: Subgroups and Covers}

 Below we recall the precise definitions of Cayley graphs
and relative Cayley graphs based on \cite{l_s}, and present Lemma
\ref{lemma1.5} from \cite{gi_sep}, which gives a characterization
of their subgraphs.

 The \emph{Cayley graph} of the group
presentation $G=gp\left\langle X|R \right\rangle$ is the oriented
graph  whose set of vertices is $G$ and whose set of edges is $G
\times X^{\pm}$, such that the edge $(g,x)$ begins at the vertex
$g$ and ends at the vertex $gx$. We denote it \fbox{$Cayley(G)$}
omitting the specification of the group presentation of $G$,
because  along this paper it is   fixed (see Section
\ref{sec:NormalForms}).

$Cayley(G)$ is a graph well-labelled  with (the alphabet)
$X^{\pm}$(that is a finite inverse automaton). Indeed, for each
edge $(g,x) \in E \left(Cayley(G) \right)$, $lab(g,x)=x$. Thus
using the automata theoretic notation, we get $g \cdot x=gx$. For
each path
$$p=( g,x_1)( g x_1,x_2)\cdots(g x_1 x_2 \cdots x_{n-1},x_n)$$ in
$Cayley(G)$, we obtain $lab(p) \equiv x_1x_2\cdots x_n \in
(X^{\pm})^*$. That is $g \cdot (x_1x_2\cdots x_n)=gx_1x_2\cdots
x_n$.

Let $S$ be a subgroup of $G=gp\left\langle X|R \right\rangle$, and
let $G/S$ denote the set of right cosets of $S$ in $G$. The
\emph{relative Cayley graph} of $G$ with respect to $S$ and the
fixed group presentation $G=gp\left\langle X|R \right\rangle$,
\fbox{$Cayley(G,S)$}
\footnote{Whenever the notation $Cayley(G,S)$ is used, it always
means that $S$ is a subgroup of the group $G$ and the presentation
of $G$ is fixed and clear from the context.}
, is an oriented graph whose vertices are the right cosets
$G/S=\{Sg \; | \; g \in G \}$, the set of edges is $(G/S) \times
X^{\pm}$, such that an edge $(Sg,x)$ begins at the vertex $Sg$ and
ends at the vertex $Sgx$.

Therefore $Cayley(G,S)$ is a graph  well-labelled with $X^{\pm}$
such that for each edge $(Sg,x)$ in $Cayley(G,S)$, $lab(Sg,x)=x$.
Using the automata theoretic notation, we get $(Sg) \cdot x=Sgx$.
Thus for each path
$$p=(S g,x_1)(S g x_1,x_2)\cdots(S g x_1 x_2 \cdots x_{n-1},x_n)$$
 in $Cayley(G,S)$, $lab(p) \equiv x_1x_2\cdots x_n \in
(X^{\pm})^*$ and $(Sg) \cdot (x_1 \cdots x_n)=Sgx_1 \cdots x_n$.

Any path $p$ in $Cayley(G,S)$ which begins at  $S \cdot 1$
\footnote{We write $S \cdot 1$ instead of the usual $S1=S$ to
distinguish this vertex of $Cayley(G,S)$ as the basepoint of the
graph.}
must end at $S~lab(p)$, so $p$ is a closed path at $S \cdot 1$ if
and only if $lab(p) \in S$. Therefore, $$Lab(Cayley(G,S),  S \cdot
1)=S.$$

$S$ acts on the Cayley graph of $G$ by left multiplication, and
$Cayley(G,S)$ can be defined as the quotient of the Cayley graph
of $G$ by this action.

Let $K$ be the standard 2-complex presenting the group
$G=\left\langle X|R \right\rangle$ (see \cite{stil}, p. 157, for
the precise definition). Thus $K$ has one vertex, $|X|$ oriented
edges and $|R|$ 2-cells. As is well known (see \cite{stil},
pp.162-163), a geometric realization of a relative Cayley graph of
$G$ is a 1-skeleton of a topological cover of $K$. This enables us
to call relative Cayley graphs of $G$,  \emph{``covers of $G$''}.

One sees, that $Cayley(G, S)$ is (the 1-skeleton of) a
finite-sheeted cover (of $K$) if and only if it has a finite
number of vertices, which means that $S$ has finite index in $G$
(\cite{stil}, p. 162). However, the generating set $X$ of $G$
might be infinite, and then a finite-sheeted cover of $G$ is an
infinite graph. Thus the term ``finite cover'' is problematic in
general. Nevertheless all groups which appear in this paper are
finitely generated. This make it possible to use the above
terminology without confusion.

The following result of Gitik \cite{gi_sep}  gives a
characterization of subgraphs of relative Cayley graphs. In order
to state it,  the definitions below are needed.

A labelled graph $\Gamma$ is \emph{$G$-based},  if any path $p$ in
$\Gamma$ with $lab(p)=_G 1_G$ is closed.
Thus any $G$-based graph is necessarily well-labelled.

Let $x \in X^{\pm}$ and $v \in V(\Gamma)$. The graph $\Gamma$ is
\emph{$x$-saturated} at $v$, if there exists $e \in E(\Gamma)$
with $\iota(e)=v$ and $lab(e)=x$.  $\Gamma$ is
\emph{$X^{\pm}$-saturated} if it is $x$-saturated for each $x \in
X^{\pm}$ at each  $v \in V(\Gamma)$.

\begin{lem}[Lemma 1.5 in \cite{gi_sep}] \label{lemma1.5}
Let $G=gp\langle X|R \rangle$ be a group and let $(\Gamma,v_0)$ be
a graph well-labelled with $X^{\pm}$. Denote $Lab(\Gamma,v_0)=S$.
Then
\begin{itemize}
    \item $\Gamma$ is $G$-based if and only if it can be embedded in $(Cayley(G,S), S \cdot 1)$,
    \item $\Gamma$ is $G$-based and $X^{\pm}$-saturated if and only if it is isomorphic to  \linebreak[4] $(Cayley(G,S), S
    \cdot~1).$
\end{itemize}
\end{lem}

%%%%%%%%%%%%%%%%%%%%%%%%%%%%%%%%%%%%%%%%%%%%%%%%%%%%%%%%%%%%%%%%%%%%%%%%%%%%%%%%%%%%%%%

\section{Normal Forms and Normal Core} \label{sec:NormalForms}
Normal words in amalgams and normal paths in the corresponding
labelled graphs are our basic tools. Below we recall their
 definitions. We define the new notion of the \emph{normal
core} of $Cayley(G,H)$. This graph is canonically associated with
the subgroup $H$ and will be constructed algorithmically in
Section \ref{subsec: TheAlgorithm}.

%-------------presentation----------------------------------------------

We start by fixing the notation. From now on whenever we refer to
the group $G$ we mean the amalgam $G=G_1 \ast_{A} G_2$, and
whenever we refer to the group  presentation of $G$ we mean the
following.
We assume that the (free) factors are given by the  finite group
presentations
\begin{align} G_1=gp\langle X_1|R_1\rangle, \ \ G_2=gp\langle
X_2|R_2\rangle \ \ {\rm such \ that} \ \ X_1^{\pm} \cap
X_2^{\pm}=\emptyset. \tag{\text{$1.a$}}
\end{align}
 $A$ is a group such that there exist two
monomorphisms
\begin{align}
\phi_1:A \rightarrow G_1 \ {\rm and } \ \phi_2:A \rightarrow G_2.
\tag{\text{$1.b$}}
\end{align}
Thus $G$ has a finite group presentation
\begin{align}
G=gp\langle X_1,X_2 | R_1, R_2, \phi_1(A)=\phi_2(A) \rangle.
\tag{\text{$1.c$}}
\end{align}

We  put $X=X_1 \cup X_2$,  $R=R_1 \cup R_2 \cup
\{\phi_1(A)=\phi_2(A)\} $. Thus $G=gp\langle X|R\rangle$.

As is well known \cite{l_s, m-k-s, serre}, the free factors embed
in $G$. It enables us to identify $A$ with its monomorphic image
in each one of the free factors. Sometimes in order to make the
context clear we'll use \fbox{$G_i \cap A$}, $i \in \{1,2\}$, to
denote the monomorphic image of $A$ in $G_i$.

%By abuse of  language, sometimes we  say that $G$ is the
%\emph{amalgam of $G_1$ and $G_2$ over (the amalgamated subgroup)
%$A$}.

%------------------------------------------------------------------------

%
%-----------------------------------------------------------------
%
\subsection*{Normal Forms}

\begin{defin}[ \cite{gi_quas, l_s, serre}]
Let $G=G_1 \ast_{A} G_2$.  We say that a word $g_1g_2 \cdots g_n
\in G$ is  \underline{in normal form} if:
\begin{enumerate}
    \item [(1)] $g_i \neq_G 1$ lies in one of the free factor of $G$,
    \item [(2)] $g_i$ and $g_{i+1}$ are in different factors of $G$,
    \item [(3)] if $n \neq 1$, then $g_i \not\in A$.
\end{enumerate}
We call the sequence $(g_1, g_2, \ldots, g_n)$ a
 \underline{normal decomposition} of the element $g \in G$, where $g=_G g_1g_2 \cdots g_n$.
\end{defin}

Any $g \in G$ has a representative in a normal form, \cite{l_s,
m-k-s, serre}.  If $g \equiv g_1g_2 \cdots g_n $ is in normal form
and $n>1$, then the Normal Form Theorem \cite{l_s} implies that $g
\neq_G 1$.

By Serre \cite{serre}, if $g$ and $h$ are two different words in
normal form with normal decompositions $(g_1, g_2, \ldots,
g_{n_1})$ and $(h_1, h_2, \ldots, h_{n_2})$, respectively, then
$g=_G h$ if and only if $n_1=n_2=n$ and there exist $a_i \in A, \
1 \leq i \leq (n-1)$, such that
$$h_1=_G g_1a_1^{-1}, \ h_j=_G a_{j-1}g_ja_j^{-1},  \ 2 \leq j \leq n-1,
\ h_n=_G a_{n-1}g_n.$$
The number $n$ is unique for a given element $g$ of $G$ and it is
called the \emph{syllable length} of $g$ (the subwords $g_i$ are
called the \emph{syllables} of $g$). We denote it by
\fbox{$length(g)$}. Notice that the number of letters in the word
$g$ is called the \emph{length} of $g$ and denoted $|g|$.

Let $p$ be a path in the graph $\Gamma$, and let $$p_1p_2 \cdots
p_n$$ be its decomposition into maximal monochromatic subpaths
(i.e., subpaths labelled with either $X_1^{\pm}$ or $X_2^{\pm}$).
Following the notation of Gitik, \cite{gi_sep}, we say that $p$ is
in \emph{normal form} (by abuse of language, $p$ is a \emph{normal
path}) if the word
$$lab(p) \equiv lab(p_1)lab(p_2) \cdots lab(p_n)$$ is in normal form.

If each $p_i$, $1 \leq i \leq n$ is a \emph{geodesic} in
$Cayley(G_j)$ (a \emph{geodesic} is the shortest path joining two
vertices) $j \in \{1,2\}$, we say that $p$ is in \emph{strong
normal form} (i.e. a \emph{strong normal path}).

%%-----------------------------------------
%
\subsection*{Normal Core}
%
%

%%%%%%%%%%%%%%%%%%%%%%%%%%%%%%%%%%%%%%%%%%%%%%%%%%%%%%%%%%%%%%%%%

\begin{defin} \label{def: normal core}
A vertex of $Cayley(G,H)$ is called \underline{essential} if there
exists a normal path closed at $H \cdot 1$ that goes through it.

The \underline{normal core} $(\Delta, H \cdot 1)$ of $Cayley(G,H)$
is the restriction of $Cayley(G,H)$ to the set of all
\underline{essential vertices}.
\end{defin}

\begin{remark} \label{rm core=union of closed paths_finite grps}
{\rm
    Note that the normal core $(\Delta, H \cdot 1)$ can be viewed as the union
    of all normal paths closed at $H \cdot 1$  in $(Cayley(G,H), H \cdot 1)$.
    Thus $(\Delta, H \cdot 1)$ is a connected graph with  basepoint $H \cdot
    1$.

    Moreover, $V(\Delta)=\{H \cdot 1\}$ and  $E(\Delta)=\emptyset$ if and only if $H$ is the trivial
    subgroup. Indeed, $H$ is not trivial iff there exists $1
    \neq g \in H$ in  normal form iff there exists $1
    \neq g \in H$ such that $g$ labels a normal path in $Cayley(G,H)$
    closed at $H \cdot 1$, iff  $E(\Delta) \neq \emptyset$.}

    \e
\end{remark}

%%%%%%%%%%%%%%%%%%%%%%%%%%%%%%%%%%%%%%%%%%%%%%%%%%%%%%%%%%%%%%%%%%%%%%%%%%%%%%%%%%%

\section{Reduced Precovers}

The notion of \emph{precovers} was defined by Gitik in
\cite{gi_sep} for subgroup graphs of amalgams. Such graphs can be
viewed as a part of the corresponding covers of $G$, that explains
the use of the term  ``precovers''. Precovers are interesting from
our point of view, because, by Lemma \ref{lemma2.12}, they allow
reading off normal forms on the graph. However these graphs could
have (\emph{redundant}) monochromatic components such that no
closed normal path starting at the basepoint goes through them.
Therefore, when looking for normal forms,  our attention can be
restricted to precovers with no redundant monochromatic components
-- \emph{reduced precovers}.

\subsection*{Precovers} \label{sec: precovers}
We say that a vertex $v \in V(\Gamma)$ is \emph{bichromatic} if
there exist edges $e_1$ and $e_2$ in $\Gamma$ with
$$\iota(e_1)=\iota(e_2)=v \ {\rm and} \ lab(e_i) \in X_i^{\pm}, \ i \in \{1,2\}.$$
The  set of bichromatic vertices of $\Gamma$ is denoted by
$VB(\Gamma)$. The vertex $v \in V(\Gamma)$ is called
\emph{$X_i$-mononochromatic} if all the edges of $\Gamma$
beginning at $v$ are labelled with $X_i^{\pm}$. We denote the set
of $X_i$-monochromatic vertices of $\Gamma$ by $VM_i(\Gamma)$ and
put $VM(\Gamma)= VM_1(\Gamma) \cup VM_2(\Gamma)$.

A subgraph of $\Gamma$ is called \emph{monochromatic} if it is
labelled only with $X_1^{\pm}$ or only with $X_2^{\pm}$. An
\emph{$X_i$-monochromatic component} of $\Gamma$ ($i \in \{1,2\}$)
is a maximal connected subgraph of $\Gamma$ labelled with
$X_i^{\pm}$, which contains at least one edge. Recall from
Section~\ref{sec: Subgroups and Covers}, that by a \emph{cover} of
a group $G$ we mean a relative Cayley graph of $G$ corresponding
to a subgroup of $G$.

%---------------------------precover------------------------------------------------

\begin{defin}[\cite{gi_sep}] \label{def: precover}
Let $G=G_1 \ast_{A} G_2$. We say that a $ G$-based graph $\Gamma$
is a \underline{precover} of $G$ if each $X_i$-monochromatic
component of $\Gamma$ is a cover of $G_i$ ($i \in \{1,2\}$).
\end{defin}

\begin{remark}
{\rm Note that by the above definition, a precover need not be a
connected graph. However along this paper we restrict our
attention only to connected precovers. Thus any time this term
 is used, we always mean that the corresponding graph
is connected.

We follow the convention that a graph $\Gamma$ with
$V(\Gamma)=\{v\}$ and $E(\Gamma)=\emptyset$ determining the
trivial subgroup (that is $Lab(\Gamma,v)=\{1\}$) is a (an empty)
precover of $G$. }\e
\end{remark}

\begin{figure}[!h]
\psfrag{x }[][]{$x$} \psfrag{y }[][]{$y$} \psfrag{v }[][]{$v$}
\psfrag{u }[][]{$u$}
\psfrag{w }[][]{$w$}
\psfrag{x1 - monochromatic vertex }[][]{{\footnotesize
$\{x\}$-monochromatic vertex}}
\psfrag{y1 - monochromatic vertex }[][]{\footnotesize
{$\{y\}$-monochromatic vertex}}
\psfrag{ bichromatic vertex }[][]{\footnotesize {bichromatic
vertex}}
\psfragscanon \psfrag{G }[][]{{\Large $\Gamma_1$}}
\psfragscanon \psfrag{K }[][]{{\Large $\Gamma_2$}}
\psfragscanon \psfrag{H }[][]{{\Large $\Gamma_3$}}
\psfragscanon \psfrag{L }[][]{{\Large $\Gamma_4$}}
\includegraphics[width=\textwidth]{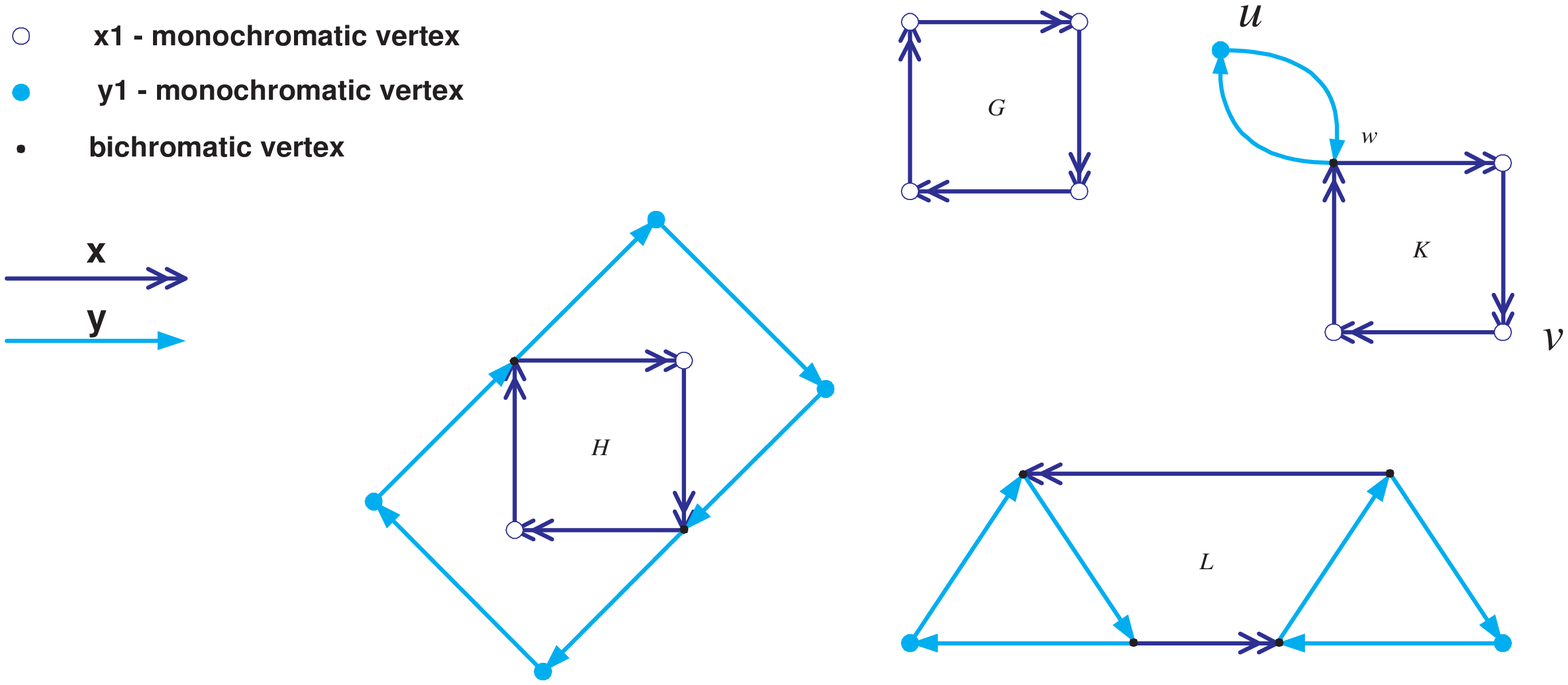}
\caption{ \label{fig:Precovers}}
\end{figure}

\begin{ex}
{\rm
 Let $G=gp\langle x,y | x^4, y^6, x^2=y^3 \rangle=\mathbb{Z}_4 \ast_{\mathbb{Z}_2} \mathbb{Z}_6$.

Recall that $G$ is isomorphic to $SL(2,\mathbb{Z})$ under the
homomorphism
$$x\mapsto \left(
\begin{array}{cc}
0 & 1 \\
-1 & 0
\end{array}
\right), \ y \mapsto \left(
\begin{array} {cc}
0 & -1\\
1 & 1
\end{array}
 \right).$$
The graphs $\Gamma_1$ and $\Gamma_3$ on Figure \ref{fig:Precovers}
are examples of precovers of $G$ with one monochromatic component
and two monochromatic components, respectively.

Though the $\{x\}$-monochromatic component of the graph $\Gamma_2$
is a cover of $\mathbb{Z}_4 $ and the $\{y\}$-monochromatic
component is a cover of $\mathbb{Z}_6$, $\Gamma_2$ is not a
precover of $G$, because it is not a $G$-based graph: $v \cdot
(x^2y^{-3})=u$, while $x^2y^{-3}=_G 1$.

The graph $\Gamma_4$ is not a precover of $G$ because its
$\{x\}$-monochromatic components are not covers of  $\mathbb{Z}_4
$. }\e
\end{ex}

\begin{remark}
{\rm Let $\Gamma$ be a precover of $G$ with $Lab(\Gamma,v_0)=H
\leq G$. By Lemma~\ref{lemma1.5}, $\Gamma$ is a subgraph of
$Cayley(G,H)$.}\e
%
%Conversely, a subgraph of $Cayley(G,H)$ is a precover of $G$. ,
%$H=Lab(\Gamma,v_0)$, is a precover if and only if each
%$X_i$-monochromatic component of $\Gamma$ is a cover of $G_i$, for
%all $i \in \{1,2\}$.}
%
%\hfill $\diamond$
%
\end{remark}

\begin{remark} \label{remark: morphism of precovers}
{\rm Let $\phi: \Gamma \rightarrow \Delta$ be a morphism of
labelled graphs. If $\Gamma$ is a precover of $G$, then
$\phi(\Gamma)$ is a precover of $G$ as well.

Indeed, a morphism of labelled graphs preserves labels and
commutes with endpoints. Thus $v \in V(\Gamma)$ is
$X_1^{\pm}$-saturated/$X_2^{\pm}$-saturated/$X_1^{\pm} \cup
X_2^{\pm}$-saturated  implies $\phi(v) \in V(\Delta)$ is
$X_1^{\pm}$-saturated/$X_2^{\pm}$-saturated/$X_1^{\pm} \cup
X_2^{\pm}$-saturated. Furthermore, let $\phi(p)$ be a path in
$\phi(\Gamma)$ with $lab(\phi(p))=_G 1$. Therefore the path $p$ in
$\Gamma$ satisfies $lab(p)=lab(\phi(p))=_G 1$. Since $\Gamma$ is a
precover, $p$ is closed. Hence the path $\phi(p)$ is closed in
$\Delta$. Therefore $\phi(\Gamma)$ is $G$-based. In particular,
$\phi(\Gamma)$ is $G_i$-based, $i \in \{1,2\}$. By Lemma
\ref{lemma1.5}, each $X_i$-monochromatic component of
$\phi(\Gamma)$ is a cover of $G_i$. Hence $\phi(\Gamma)$ is a
precover of $G$. }\e
\end{remark}

%--------------------------------compatibility------------------------------------

The graph $\Gamma$  is called \emph{compatible at a bichromatic
vertex} $v$ if for any monochromatic path $p$ in $\Gamma$ such
that $\iota(p)=v$ and $lab(p) \in A$ there exists a monochromatic
path $t$ of a different color in $\Gamma$ such that $\iota(t)=v$,
$\tau(t)=\tau(p)$ and $lab(t)=_G lab(p)$. We say that $\Gamma$ is
\emph{compatible} if it is compatible at all bichromatic vertices.

\begin{ex}
{\rm The graphs $\Gamma_1$ and $\Gamma_3$ on Figure
\ref{fig:Precovers} are compatible. The graph $\Gamma_2$ does not
possess this property because $w \cdot x^{2}=v$, while $w \cdot
y^3=u$. $\Gamma_4$ is not compatible as well.} \e
\end{ex}

\begin{remark} [Remark 2.11 in \cite{gi_sep}] \label{remark:
precovers are compatible} {\rm Precovers are compatible. \hfill
$\diamond$}
\end{remark}

\begin{remark} [Corollary2.13 in \cite{gi_sep}]  \label{corol2.13}
Let $\Gamma$ be a compatible graph. If all  $X_i$-components of
$\Gamma$ are $G_i$-based, $i \in \{1,2\}$, then $\Gamma$ is
$G$-based. In particular, if each $X_i$-component of $\Gamma$ is a
cover of $G_i$, $i \in \{1,2\}$, and $\Gamma$ is compatible, then
$\Gamma$ is a precover of $G$.
\end{remark}

%%%%%%%%%%%%%%%%%%%%%%%%%%%%%%%%%%%%%%%%%%%%%%%%%%%%%%%%%%%%%%%%%%%%%
%
% \section{Reduced Precovers} \label{subsec:ReducedPrecovers}
%
%-------------------------------------------------------------------

%%%%%%%%%%%%%%%%%
% Summarizing the above discussion,  one can think about a precover
% of $G$ as a bunch of factors covers  which are glued to each other
% according to the amalgamation (such that the resulting graph is
% connected, well-labelled with $X_1^{\pm} \cup X_2^{\pm}$ and
% compatible with $A$).
%%%%%%%%%%%%%%%%%

Recall that our objective is to be able to read normal words on
the constructed graph. The following lemma of Gitik shows that
precovers are suitable for this purpose.

\begin{lem} [Lemma2.12 in \cite{gi_sep}] \label{lemma2.12}
If $\Gamma$ is a compatible graph, then for any  path $p$ in
$\Gamma$ there exists a path $t$ in normal form which has the same
endpoints and the same label (in $G$) as $p$.
\end{lem}

The statement of this lemma can be even extended when the graph
$\Gamma$ is a precover.

\begin{lem} \label{lem: all normal paths are in precover}
Let $\Gamma$ be a precover of $G$. Let $p$ be a path in $\Gamma$
with $\iota(p)=v_1$, $\tau(p)=v_2$ and $lab(p) \equiv w$.

Then for each normal word $w'$ of syllable length greater than 1
such that $w'=_G w$ there exist a normal path $p'$ in $\Gamma$
with $\iota(p')=v_1$, $\tau(p')=v_2$ and $lab(p') \equiv w'$.
\end{lem}
\begin{proof}
By Lemma \ref{lemma2.12}, we can assume that the path $p$  and the
word $w$ are normal. Let $p=p_1 \cdots p_k$ be a decomposition of
$p$ into maximal monochromatic paths ($k>1$). Let $C_i$ be the
monochromatic component of $\Gamma$ containing the subpath $p_i$
of $p$ ($1 \leq i \leq k$), that is $p_i \subseteq C_i \cap p$.

Let $lab(p_i) \equiv w_i$ ($1 \leq i \leq k$). Hence $w \equiv w_1
\cdots w_k$, where $(w_1, \ldots ,w_k)$ is a normal (Serre)
decomposition of $w$ and $w_i \in G_{l_i}$ ($l_i \in \{1,2\}$).

Let $w' \equiv w'_1 \cdots w'_m$ be a normal word with the normal
(Serre) decomposition $( w'_1, \ldots, w'_m)$ such that $w=_Gw'$.
Therefore, by \cite{serre} p.4, $m=k$ and
$$w'_1=_Gw_1a_{1_{1}}^{-1}, \; \
w'_i=_Ga_{{(i-1)}_i}w_ia_{i_i}^{-1} \ (2 \leq i \leq k-1), \ \;
w'_k=_Ga_{{(k-1)}_k}w_k,$$ where  $a_{1_j}, a_{i_j}, a_{(k-1)_j}
\in A \cap G_{l_j}$ ($1 \leq j \leq k$) such that $a_{i_i}=_G
a_{i_{(i+1)}}$.

Let $u_i=\tau(p_i)$ ($1 \leq i \leq k$). Thus
 $u_i \in VB(C_i) \cap
VB(C_{i+1})$. Since $\Gamma$ is a precover of $G$, $C_i$ and
$C_{i+1}$ are covers of $G_{l_i}$ and of $G_{l_{i+1}}$,
respectively. That is they are $X_{l_i}^{\pm}$-saturated and
$X_{l_{i+1}}^{\pm}$-saturated, respectively. Hence there are paths
$t_i $ in $C_i$ and $s_{i+1}$ in $C_{i+1}$ starting at $u_i$ and
labelled by $a_{i_i}^{-1}$ and $a_{i_{i+1}}^{-1}$, respectively
(see Figure \ref{fig:AllNormalPathsInPrecovers}). Since $\Gamma$
is compatible (as a precover of $G$), $\tau(t_i)=\tau(s_{i+1})$.
%
%---------------------------------------
\begin{figure}[!h]
\psfrag{p1 }{$p_1$} \psfrag{pi }{$p_{i}$}  \psfrag{pk }{$p_k$}
\psfrag{q1 }[][]{$p'_1$} \psfrag{qi }[][]{$p'_i$} \psfrag{qk
}{$p'_k$}
\psfrag{v1 }{$v_1$} \psfrag{v2 }{$v_2$} \psfrag{u1 }{$u_1$}
\psfrag{c1 }{$C_1$} \psfrag{ci }{$C_{i}$}  \psfrag{ck }{$C_k$}
\psfrag{c2 }{$C_2$} \psfrag{cm }{$C_{k-1}$}
\psfrag{t1 }{$t_1$} \psfrag{ti }{$t_{i}$}  \psfrag{tj }{$t_{i-1}$}
\psfrag{tm }{$t_{k-1}$}
\psfrag{s2 }{$s_2$} \psfrag{si }{$s_{i}$}  \psfrag{sj }{$s_{i+1}$}
\psfrag{ sk }{$s_k$}
\includegraphics[width=\textwidth]{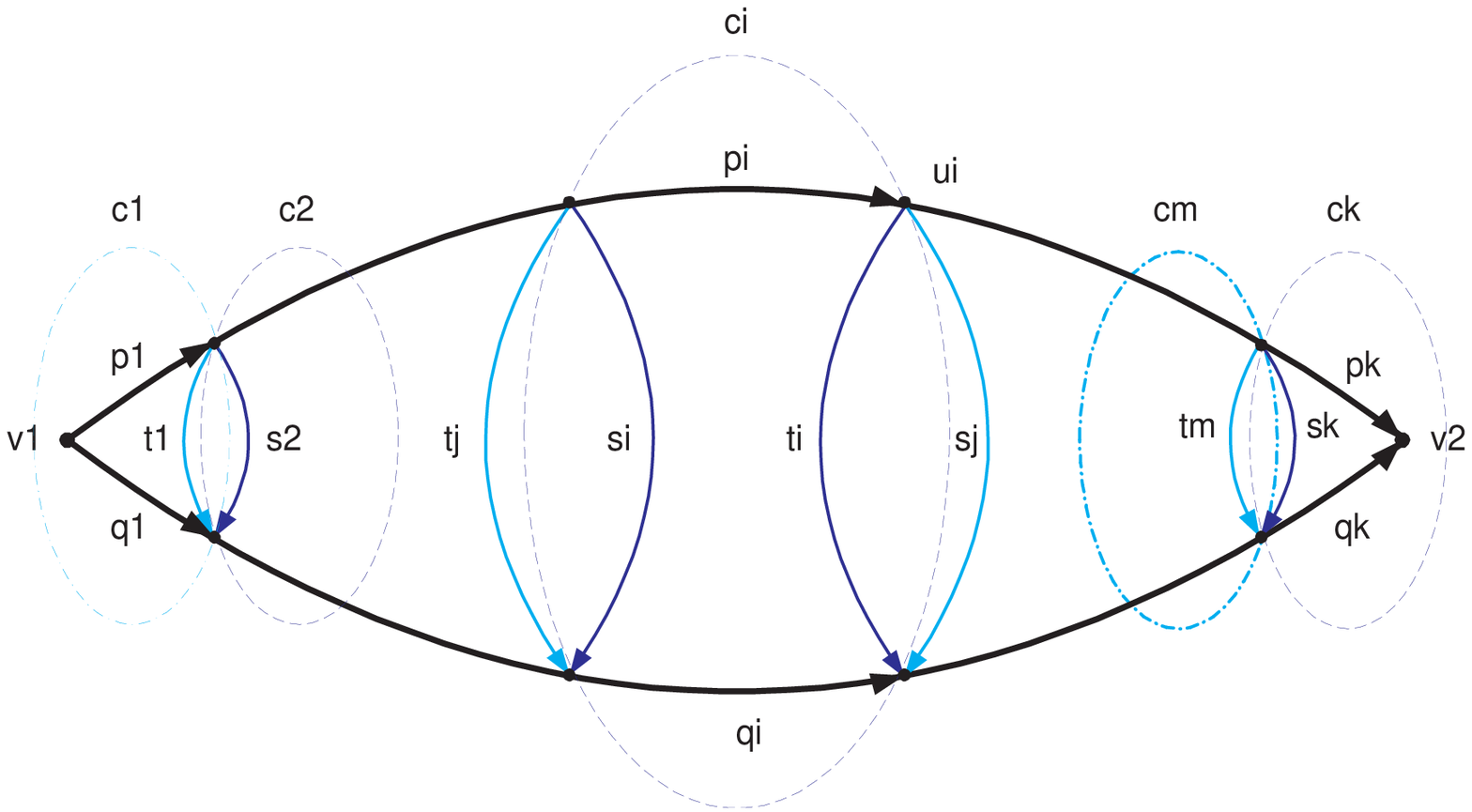}
\caption{ \label{fig:AllNormalPathsInPrecovers}}
\end{figure}
%---------------------------------------
%
Hence there exists a path $\gamma$ in $\Gamma$ such that
$\gamma=\gamma_1 \cdots \gamma_k$, where
$$\gamma_1=p_1{t_1} \subseteq C_1, \ \; \gamma_i=\overline{s_i}p_i{t_i} \subseteq C_i \ (2 \leq i \leq k-1),
\ \; \gamma_k=\overline{s_k}p_k \subseteq C_k.$$
Thus $\iota(\gamma)=v_1$, $\tau(\gamma)=v_2$ and
$$lab(\gamma_1) \equiv w_1a_{1_1}^{-1}, \ \; lab(\gamma_i) \equiv
a_{(i-1)_i}w_ia_{i_i}^{-1}  \ (2 \leq i \leq k-1), \ \;
lab(\gamma_k) \equiv
a_{(k-1)_k}w_k.$$ % Hence $lab(\gamma)=_G w$.
Since $w'_i=_G lab(\gamma_i)$ ($1 \leq i \leq k$) and because the
component $C_i$ is $X_{l_i}^{\pm}$-saturated, there exists a path
$p'_i$ in $C_i$ such that $\iota(p'_i)=\iota(\gamma_i)$ and
$lab(p'_i)\equiv w'_i$. Moreover, $\tau(p'_i)=\tau(\gamma_i)$,
because the component $C_i$ is $G_{l_i}$-based. Therefore there
exists a path $p'=p'_1 \cdots p'_k$ in $\Gamma$ such that
$\iota(p')=v_1$, $\tau(p')=v_2$ and $lab(p') \equiv w'$.

\end{proof}

\begin{figure}[!htb]
\begin{center}
\psfrag{v0 }[][]{$v_0$}
\psfrag{x }[][]{$x$} \psfrag{y }[][]{$y$}
\psfragscanon \psfrag{G }[][]{{\Large $\Gamma$}}
\psfrag{H }[][]{\Large $\Gamma'$} \psfrag{C }[][]{\Large $C$}
\psfrag{A }[][]{\Large $D$}
\includegraphics[width=0.7\textwidth]{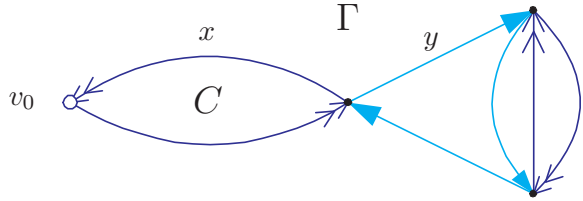}
\caption{ {\footnotesize We use the same labelling as on
Figure~\ref{fig:Precovers}} \label{fig:NormalPathsPrecovers}}
\end{center}
\end{figure}

%--------------------------------------------------------------------
\begin{remark} \label{remark:NormalPathsPrecovers}
{\rm When $length(w)=1$ the statement of Lemma~\ref{lem: all
normal paths are in precover} need not be true. Thus, for example,
the graph $\Gamma$, illustrated on
Figure~\ref{fig:NormalPathsPrecovers}, is a precover of
$G=gp\langle x,y | x^4, y^6, x^2=y^3 \rangle=\mathbb{Z}_4
\ast_{\mathbb{Z}_2} \mathbb{Z}_6$. There is a path $p$ in $\Gamma$
with $lab(p) \equiv x^2$ and $\iota(p)=\tau(p)=v_0$. However there
is no path $p'$ in $\Gamma$ with the same endpoints as $p$ and
$lab(p') \equiv y^3$. } \e
\end{remark}

%---------------------------------------------------------------------

\begin{cor} \label{cor: distance between normal paths}
Let $p$ and $p'$ be as in Lemma~\ref{lem: all normal paths are in
precover}.

If $G=G_1 \ast_A G_2$ is an amalgam of finite groups $G_1$ and
$G_2$ then $p \subset N_d(p')$ and $p' \subset N_d(p)$, where
$d=max(diameter(G_1), diameter(G_2))$.
\end{cor}
\begin{proof}

Recall that a {\it group diameter} is the length of the longest
geodesic in its Cayley graph.

Thus $d_j=diameter(G_j)=diameter(Cayley(G_j))$ ($j=1,2$) is
finite.
Since each $X_{j}$-monochromatic component $C$ of $\Gamma$ is
isomorphic to either $Cayley(G_j)$, $diameter(C)=d_j$. Let
$d=max(d_1, d_2)$.

By the proof of Lemma \ref{lem: all normal paths are in precover},
$p_i \subseteq C_i$ and $p_i' \subseteq C_i$. Thus $p_i \subset
N_{d}(p'_i)$ and $p'_i \subset N_{d}(p_i)$ ($ 1 \leq i \leq k$).
Hence $p \subset N_{d}(p')$ and  $p' \subset N_{d}(p)$.
\end{proof}

However some monochromatic components of precovers may carry no
``essential information''  concerning normal forms. More
precisely, if in a monochromatic component $C$ of the precover
$\Gamma$ every path between any two bichromatic vertices of $C$ is
labelled by an element of $A$, then, evidently, no normal path in
$\Gamma$ goes through this component, see Figure
\ref{fig:RedundantComponent} and Example
\ref{ex:RedundantComponent}.

Below we are looking for an explicit characterization of such
(\emph{redundant}) monochromatic components of precovers. This
enables us to present the new notion of a \emph{reduced precover},
which is, roughly speaking, a precover with no redundant
monochromatic components.

%----------------------------------------------------------------
%
\subsection*{Redundant Monochromatic Components}
Let $\Gamma$ be a precover of $G$.  Let $C$ be a
$X_i$-monochromatic component of $\Gamma$ ($i \in \{1,2\}$). Then
$A$ \emph{acts} on $V(C)$ by right multiplication.

Let $v \in V(C)$, then the \emph{$A$-orbit} of $v$ is
$$A(v)=\{v \cdot a \; | \; a \in A\}.$$
Since $\Gamma$ is a precover of $G$, it is compatible with $A$.
Thus $v \in VB(C)$ if and only if $A(v)\subseteq VB(C)$. Hence
bichromatic vertices of $C$ are determined by the
\emph{$A$-action}. Moreover, $A(v)=VB(C)$ if and only if the
bichromatic vertices of $C$ form the unique $A$-orbit.

\begin{claim} \label{claim: red precover}
For all $v_1,v_2 \in VB(C)$, $v_1 \cdot a=v_2$ implies $a \in A$
if and only if $VB(C)=A(\vartheta)$ and $Lab(C,\vartheta)=K \leq
A$, for all $\vartheta \in VB(C)$.
\end{claim}

In other words, each path $p$ in $C$ ($ C \subseteq \Gamma$) with
$\iota(p), \tau(p) \in VB(C)$ satisfies  $lab(p) \in A$ if and
only if there exists a unique $A$-orbit of bichromatic vertices in
$C$ and $Lab(C,\vartheta) \leq A$, for all $\vartheta \in VB(C)$.

\begin{figure}[!htb]
\begin{center}
\psfragscanon \psfrag{A }[][]{{\Large $\Gamma_1$}}
\psfrag{v }[][]{$v$} \psfrag{u }[][]{$u$}
\psfrag{B }[][]{\Large $\Gamma_2$}
\psfrag{G }[][]{\Large $\Gamma_3$}
\psfrag{H }[][]{\Large $\Gamma_4$}
 \psfrag{C }[][]{\large $C$} \psfrag{D }[][]{\large $D$}
\includegraphics[width=0.9\textwidth]{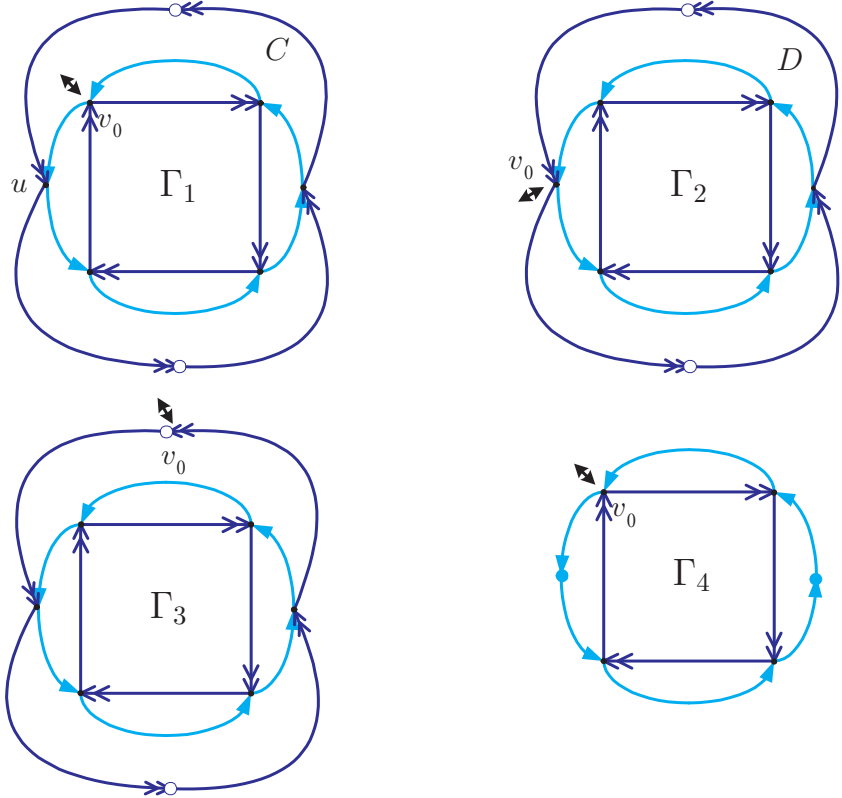}
\caption{ {\footnotesize We use the same labelling as on
Figure~\ref{fig:Precovers}} \label{fig:RedundantComponent}}
\end{center}
\end{figure}

\begin{proof}[Proof of Claim \ref{claim: red precover}]
Assume first that $VB(C)=A(\vartheta)$ and $K=Lab(C,\vartheta)
\leq A$. Let $v_1 , v_2 \in VB(C)$. Since $(C,\vartheta)$ is
isomorphic to $Cayley(G_i, K, K \cdot 1)$ and $C$ has the unique
$A$-orbit of bichromatic vertices, there exist $a_1, a_2 \in A$
such that $v_1=(K \cdot 1) \cdot a_1=Ka_1$ and $v_2=(K \cdot 1)
\cdot a_2=Ka_2$. Thus
$$ v_1 \cdot a= v_2 \ \Leftrightarrow \ (K a_1) \cdot a= K a_2 \ \Leftrightarrow \
a_1aa_2^{-1} \in K.$$ Since $K \leq A$, we have $a \in A$.

Conversely, assume that for each pair of vertices $v_1,v_2 \in
VB(C)$ each path $p$ in $C$ with $\iota(p)=v_1$ and $\tau(p)=v_2$
has $lab(p) \equiv a \in A$. In particular, if $v_1=v_2=\vartheta
\in VB(C)$ then $\vartheta \cdot x =\vartheta$ implies $x \in A$.
However  $ x \in Lab(C,\vartheta)=K$. Therefore
$Lab(C,\vartheta)=K \leq A$. The equality $VB(C) = A(\vartheta)$
holds by the definition of $A$-orbits, because $\vartheta \in
VB(C)$.

\end{proof}
 Now we are ready to give a precise definition of the new notion of \emph{redundant monochromatic components}.
\begin{defin} \label{def: redundant component}
Let $(\Gamma,v_0)$ be a precover of $G$.
Let $C$ be a $X_i$-monochromatic component of $\Gamma$ ($i \in
\{1,2\}$).  $C$ is  \underline{redundant} if one of the following
holds.
\begin{enumerate}
\item [(1)] $C$ is the unique monochromatic component of $\Gamma$
(that is $\Gamma=C$) and $Lab(C,v_0)=\{1\}$ (equivalently, by
Lemma \ref{lemma1.5}, $C$ is isomorphic to $Cayley(G_i)$).
\item [(2)] $\Gamma$ has at least two distinct monochromatic
components and the following holds.

Let $\vartheta \in VB(C)$. Let $K=Lab(C,\vartheta)$ (equivalently,
by Lemma \ref{lemma1.5}, $(C,\vartheta)=(Cayley(G_i,K), K \cdot
1)$). Then
\begin{itemize}
    \item [(i)] $K \leq A$,
    \item [(ii)] $VB(C)=A(\vartheta)$,
%    \item [(iii)] either $v_0 \not\in VM(C)$,  $K=\{1\}$ or, $v_0 \not\in V(C)$ and $K \neq \{1\}$.
\item [(iii)] either $v_0 \not\in V(C)$ or, $v_0 \in VB(C)$ and $K
= \{1\}$.
\end{itemize}
\end{enumerate}
\end{defin}
\begin{ex} \label{ex:RedundantComponent}
{\rm Let $G=gp\langle x,y | x^4, y^6, x^2=y^3 \rangle=\mathbb{Z}_4
\ast_{\mathbb{Z}_2} \mathbb{Z}_6$.

The graphs on Figure \ref{fig:RedundantComponent} are examples of
precovers of $G$. The $\{x\}$-monochromatic component $C$ of the
graph $\Gamma_1$ is redundant, because $(C,u)$ is isomorphic to
$Cayley(\mathbb{Z}_4)$, that is $Lab(C,u)=\{1\}$, while
$|VB(C)|=2=[\mathbb{Z}_4 : \mathbb{Z}_2]$ and $v_0 \not\in V(C)$.

The $\{x\}$-monochromatic component $D$ of the graph $\Gamma_2$ is
redundant, because $Lab(D,v_0)=\{1\}$, while $v_0 \in
VM(\Gamma_2)$.

However the graphs $\Gamma_3$ and $\Gamma_4$ have no redundant
components. } \e
\end{ex}
%

%-------------------------------------------------

\begin{remark} \label{remark: comput cond for redundant
component}
{\rm
Note that when the free factors $G_1$ and $G_2$ of the amalgam
$G=G_1 \ast_A G_2$ are finite groups, then Claim \ref{claim: red
precover} and Definition \ref{def: redundant component} can be
restated in the following computational manner.

Recall that the \emph{$A$-stabilizer} of $v$ is
$$A_v=\{ a \in A \; | \; v \cdot a =v \} \leq A.$$
As is well-known, the cosets of the stabilizer subgroup are in a
one-to-one correspondence with the elements in the orbit
$$A(\vartheta) \sim A/A_{\vartheta}.$$

Thus
$$|A(\vartheta)|=[A:A_{\vartheta}].$$
Let $\vartheta \in VB(C)$. Let $K=Lab(C,\vartheta)$ (equivalently,
$(C,\vartheta)=(Cayley(G_i,K), K \cdot 1)$). Hence  $A_{\vartheta}
= K \cap A$.
Since $VB(C) = A(\vartheta)$ if and only if
$|VB(C)|=|A(\vartheta)|$, the condition  $K \leq A$ implies $VB(C)
= A(\vartheta)$ if and only if $|VB(C)|=[A:K]$.

This enables us to replace the condition $VB(C) = A(\vartheta)$ in
Claim \ref{claim: red precover} and in Definition \ref{def:
redundant component} by its computational analogue
$|VB(C)|=[A:K]$.

}\e
\end{remark}

%
%

%-------------------------------------------------------------------

Let us make the following assumption. From now on whenever we say
that a path $p$ in $\Gamma$ \emph{goes through the vertex} $v \in
V(\Gamma)$, we mean that $v \in V(p)$. And whenever we say that a
path $p$ in $\Gamma$  \emph{goes through the monochromatic
component} $C$ in $\Gamma$, we mean that $E(p) \cap E(C) \neq
\emptyset$. That is if $p=p_1 \cdots p_n$ is a decomposition of
$p$ into maximal monochromatic paths then there exists $1 \leq l
\leq n$ such that $C$ contains the subpath $p_l$ ($p_l \subseteq p
\cap C$ or, more precisely, $E(p_l) \subseteq E(p) \cap E(C)$).

\begin{lem}  \label{red-precover=normal paths}
Let $(\Gamma,v_0)$ be a precover of $G$. Then $X_i$-monochromatic
component $C$ of $\Gamma$ ($i \in \{1,2\}$) is redundant if and
only if no normal path $p$ in $\Gamma$ closed at $v_0$ goes
through $C$.
\end{lem}
\begin{proof}
Let $C$ be a $X_i$-monochromatic component of $\Gamma$ ($i \in
\{1,2\}$).
Let $p$ be a path closed at $v_0$ that goes through $C$. Let
$p=p_1p_2 \cdots p_k$ be a decomposition of $p$ into maximal
monochromatic paths. Thus there exists $1 \leq j \leq k$ such that
$p_j \subseteq C$.

If $k=1$ then  $p \subseteq C$ and $v_0 \in V(C)$. Thus $p$ is
normal if and only if $lab(p) \neq \{1\}$ if and only if
$Lab(C,v_0) \neq \{1\}$ if and only if neither condition (1) nor
condition (iii) in  Definition \ref{def: redundant component} is
satisfied.

Assume now that $k > 1$. The path $p$ is normal if and only if
$lab(p_j) \not\in A$ for all $1 \leq j \leq k$.
%
%
%If $p_1 \subseteq C$ or $p_k \subseteq C$ then $lab(p_j) \not\in
%A$ ($j \in \{1,k\}$)  if and only if $v_0 \not\in VB(C)$, or $v_0
%\in VB(C)$ and, by Claim \ref{claim: red precover}, conditions
%$(i)$ and $(ii)$ in Definition \ref{def: redundant component} are
%not satisfied. If $p_j \subseteq C$ ($2 \leq j \leq k-1$) then
%$lab(p_j) \not\in A$ if and only if conditions $(i)$ and $(ii)$ in
%Definition \ref{def: redundant component} are not satisfied, by
%Claim \ref{claim: red precover}.
%
By Claim \ref{claim: red precover}, this happens if and only if at
least one of the conditions $(i)$, $(ii)$ in Definition \ref{def:
redundant component} is not satisfied for the monochromatic
component $C$.

Therefore $p$ is normal if and only if $C$ is not redundant.

\end{proof}

%-----------------------------------------------------------------------

Now we show that removing of a redundant monochromatic component
from a precover $(\Gamma,v_0)$ leaves the resulting graph  a
precover and don't change the subgroup determined by the graph.

One can think of this procedure as an analogue of the   ``cutting
hairs'' procedure, presented by Stallings in \cite{stal}, for
subgroup graphs in the case of free groups. Indeed, a hair is cut
from the graph because no freely reduced paths closed at the
basepoint go through the hair. Similarly, when interested in
normal paths closed at the basepoint of a precover, its redundant
components can be erased, because no such paths go through them.

\begin{figure}[!htb]
\begin{center}
\psfragscanon \psfrag{G }[][]{{\Large $\Gamma$}}
\psfrag{H }[][]{\Large $\Gamma'$}

\psfrag{v0 }[][]{$v_0$}
\includegraphics[width=0.9\textwidth]{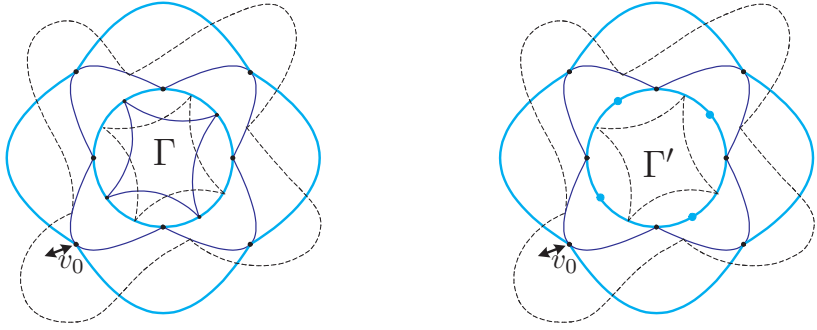}
\caption{{\footnotesize The closed grey curves represent
$G_1$-monochromatic components. The closed black curves represent
 $G_2$-monochromatic components. The broken lines denote the
rest of the graphs $\Gamma$ and $\Gamma'$, respectively. The small
black dots are bichromatic vertices. The grey dots are
$X_1$-monochromatic vertices.}
\label{fig:RemovingRedundantComponent}}
\end{center}
\end{figure}

Let $(\Gamma, v_0)$ be a precover of $G$.  Let $C$ be a redundant
$X_j$-monochromatic component of $\Gamma$ ($j \in \{1,2\}$). We
say that the graph $\Gamma'$ is obtained from the graph $\Gamma$
by \emph{removing of redundant $X_j$-monochromatic component} $C$
, if $\Gamma'$ is obtained by removing all edges and all
$X_j$-monochromatic vertices of $C$, while keeping all its
bichromatic vertices (see Figure
\ref{fig:RemovingRedundantComponent}). More precisely, if
$\Gamma=C$ then we set $V(\Gamma')=\{v_0\}$,
$E(\Gamma')=\emptyset$. Otherwise $V(\Gamma') = V(\Gamma)
\setminus VM_j(C)$, where
\begin{eqnarray}
VB(\Gamma') & = & VB(\Gamma) \setminus VB(C),  \nonumber \\
VM_j(\Gamma') & = & VM_j(\Gamma) \setminus VM_j(C), \ (1 \leq i
\neq j
\leq 2), \nonumber \\
 VM_i(\Gamma') & = & VM_i(\Gamma) \cup VB(C).  \nonumber
\end{eqnarray}
And $$E(\Gamma')=E(\Gamma) \setminus E(C) \ {\rm and} \
lab_{\Gamma'}(e) \equiv lab_{\Gamma}(e) \; (\forall e \in
E(\Gamma')).  $$

%------------------------------------------------------
%

%-----------------------------------------------------

\begin{lem} \label{Cayley-remove}
 Let $(\Gamma, v_0)$ be a precover of $G$.
Let $\Gamma'$ be the graph obtained from the graph $\Gamma$ by
removing of a redundant $X_j$-monochromatic component $C$ of
$\Gamma$ ($j \in \{1,2\}$).

Then $\Gamma'$ is a precover of $G$ such that
$Lab(\Gamma,v_0)=Lab(\Gamma',v_0')$, where $v_0$ is the basepoint
of $\Gamma$ and $v_0'$ is the (corresponding) basepoint of
$\Gamma'$.
\end{lem}
\begin{proof}
If $\Gamma=C$ then $\Gamma'$ is a precover, by the construction.
Since $C$ is redundant, $Lab(\Gamma,v_0)=\{1\}$. On the other
hand, $Lab(\Gamma',v_0)=\{1\}$ as well.  We are done.

Assume now that $\Gamma$ has at least two monochromatic
components.

Evidently, by the construction, $\Gamma'$ is a precover of $G$.
Indeed, since $VB(\Gamma') \subseteq VB(\Gamma)$ and $\Gamma$ is
compatible, then so is  $\Gamma'$. Let $D \neq C$ be a
monochromatic component of $\Gamma$. Then $D \subseteq \Gamma'$.
Thus each $X_i$-monochromatic component of $\Gamma'$ is a cover of
$G_i$ ($i \in \{1,2\}$). Hence $\Gamma'$ is a precover of $G$, by
Lemma~\ref{corol2.13}.
  Note that $(\Gamma',v'_0) \subseteq (\Gamma,v_0)$.
Thus $Loop(\Gamma',v'_0) \subseteq Loop(\Gamma,v_0)$ and we get
$Lab(\Gamma', v'_0) \subseteq Lab(\Gamma,v_0)$.

Let  $w \in Lab(\Gamma,v_0)$. Then there is $t \in
Loop(\Gamma,v_0)$ such that $lab(t)=_G w$. If no subpath of $t$ is
in $C$ then $t$ is also a path in $\Gamma'$. Therefore $lab(t)=_G
w \in Lab(\Gamma',v'_0)$.

Otherwise,  there is a decomposition $t=t_1q_1t_2q_2 \ldots
q_{k-1}t_k$ such  that $\iota(t_1)=\tau(t_k)=v_0$ and for all $1
\leq i \leq k$, $q_i$ is a path in the component $C$ and $t_i$ is
a path in $\Gamma'$ with the normal decomposition $t_i=t_{i1}
\cdots t_{im_i}$. Since $E(t_i) \cap E(C)=\emptyset$, the paths
$t_{im_i}$, $q_i$ and $q_i$, $t_{(i+1)1}$ are pairs of
monochromatic paths of different colors. Thus the vertices
$\tau(t_i)=\iota(q_i)$ and $\tau(q_i)=\iota(t_{i+1})$ are
bichromatic vertices of $\Gamma$. Therefore, since $C$ is
redundant, Claim \ref{claim: red precover} implies that $lab(q_i)
\in G_j \cap A$.

Let $D$ be a $X_l$-monochromatic component of $\Gamma$ such that
$t_{im_i}$ is a path in $ D$, where $1 \leq j \neq l \leq 2$.
Since $\Gamma$ is a precover, $D$ is a cover of $G_l$. Since the
vertex $\iota(q_i)$ is bichromatic in $\Gamma$, while $\Gamma$  is
compatible and $lab(q_i) \in G_j \cap A$, there exists a path
$p_i$ in $D$ such that
$$\iota(p_i)=\iota(q_i), \ \tau(p_i)=\tau(q_i) \ {\rm and} \
lab(p_i)=_Glab(q_i).$$

 Thus the path $ t'=t_1p_1t_2p_2
\ldots p_{k-1}t_k$ is a closed path at $v_0'$ in $\Gamma'$ with
$lab(t')=_G lab(t)$. Therefore $w \equiv lab(t)=_G lab(t') \in
Lab(\Gamma', v'_0)$.

Hence $Lab(\Gamma)=Lab(\Gamma')$.

Proceeding in the same manner as in the construction of $t'$ (in
$\Gamma'$) from the  path $t$ (in $\Gamma$), one can shows that
any two vertices of $\Gamma$ remain connected by a path in
$\Gamma'$. More precisely, given a pair of vertices $v$ and $w$ in
$\Gamma$ and given  a path $s$ in $\Gamma$ connecting them, one
can construct an appropriate path $s'$ in $\Gamma'$  such that
$\iota(s')=\iota(s)=v$, $\tau(s')=\tau(s)=w$ and $lab(s')=_G
lab(s)$. Therefore the graph $\Gamma'$ is connected.

\end{proof}

%------------------------------------
%
\subsection*{Reduced Precovers} \label{sec: reduced precovers}
\begin{defin} \label{def:reduced precover}
A precover $(\Gamma,v_0)$ of $G$ is called \underline{reduced} if
and only if the following holds
\begin{itemize}
\item [(1)] $(\Gamma,v_0)$ has no redundant monochromatic
components.
\item [(2)] If there exists a $X_i$-monochromatic component $C$ of
$\Gamma$ ($i \in \{1,2\}$) such that
 $ v_0 \in V(C) \ {\rm and } \ K \cap A \neq \{1\}, \ {\rm where}
 \ K=Lab(C,v_0)$ (equivalently,
$(C,v_0)=(Cayley(G_i,K), K \cdot 1)$), then there exists  a
$X_j$-monochromatic component $D$ of $\Gamma$ ($ 1 \leq i \neq j
\leq 2$) such that $ v_0 \in V(D) \ {\rm and } \ K \cap A =_G L
\cap A, \ {\rm where}
 \ L=Lab(D,v_0)$ (equivalently,
$(D,v_0)=(Cayley(G_i,L), L \cdot 1)$).
\end{itemize}
\end{defin}

\begin{figure}[!htb]
\begin{center}
\psfrag{v0 }[][]{$v_0$}
\psfrag{x }[][]{$x$} \psfrag{y }[][]{$y$}
\psfragscanon \psfrag{G }[][]{{\Large $\Gamma$}}
\psfrag{H }[][]{\Large $\Gamma'$} \psfrag{C }[][]{\Large $C$}
\psfrag{A }[][]{\Large $D$}
\includegraphics[width=0.7\textwidth]{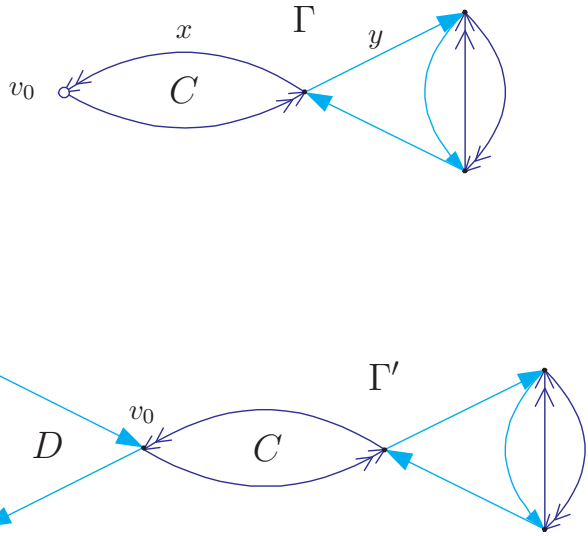}
\caption{ {\footnotesize We use the same labelling as on
Figure~\ref{fig:Precovers}}  \label{fig: ReducedPrecover}}
\end{center}
\end{figure}

\begin{remark} \label{rem: reduced precover}
{\rm Note that condition (2) in the above definition merely says
that if $A \cap H \neq \{1\}$ then $v_0 \in VB(\Gamma)$, where
$H=Lab(\Gamma,v_0)$.

Therefore if $\Gamma$ has the unique $X_i$-monochromatic component
$C$ (that is $\Gamma=C$, $i \in \{1,2\}$) then $H$ is a nontrivial
subgroup of $G_i$ such that $A \cap H = \{1\}$.

If $V(\Gamma)=\{v_0\}$ and $E(\Gamma)=\emptyset$ then $\Gamma$ is
a reduced precover, by the above definition, with
$Lab(\Gamma,v_0)=\{1\}$ } \e
\end{remark}
\begin{ex} \label{ex: reduced precover}
{\rm
 Let $G=gp\langle x,y | x^4, y^6, x^2=y^3 \rangle=\mathbb{Z}_4 \ast_{\mathbb{Z}_2} \mathbb{Z}_6$.

The precovers  $\Gamma_1$ and $\Gamma_2$ from Figure
\ref{fig:RedundantComponent} are not reduced because they have
redundant components $C$ and $D$, respectively (see Example
\ref{ex:RedundantComponent}). The graphs $\Gamma_3$ and $\Gamma_4$
from the same figure are reduced precover of $G$ because they are
precovers with no redundant components and with a bichromatic
basepoint.

The precover $\Gamma$ on Figure \ref{fig: ReducedPrecover} is not
a reduced precover of $G$ though it has no redundant components.
The problem now is the $\{x\}$-monochromatic component $C$ of
$\Gamma$ because $Lab(C,v_0)=\langle x^2 \rangle$, while the
basepoint $v_0$ is a $\{x\}$-monochromatic vertex.   It is easy to
see that the graph $\Gamma'$  obtained from $\Gamma$ by gluing at
$v_0$ the appropriate $\{y\}$-monochromatic component $D$ with
$Lab(D,v_0)=\langle y^3 \rangle$ is a reduced precover of $G$, by
Definition \ref{def:reduced precover}. } \e
\end{ex}
%
%
%---------------------------------------------------------------

Let $(\Gamma,v_0)$ be a  precover of $G$ with no redundant
components, which is not a reduced precover.  Hence $v_0 \in
VM_l(\Gamma)$ ($l \in \{1,2\}$) and the assumption of condition
(2) in Definition \ref{def:reduced precover} holds, that is
$\Gamma$ has a $X_l$-monochromatic component $C$  with
$Lab(C,v_0)=K$ such that $L=K \cap A$ is a nontrivial subgroup of
$A$.

Thus $(\Gamma,v_0)$ can be ``reconstructed'' in the obvious way
(see Figure \ref{fig: ReducedPrecover} and Example \ref{ex:
reduced precover}) such that the resulting graph is a reduced
precover of $G$ determining the same subgroup as the graph
$(\Gamma,v_0)$ does.

Let $(\Gamma',v_0')$ be the graph obtained by taking a disjoint
union of the graphs $(\Gamma,v_0)$ and $(Cayley(G_j,L),L \cdot 1)$
($1 \leq j \neq l \leq 2$) via the identification of $L \cdot 1$
with $v_0$ and the identification of the $X_j$-monochromatic
vertices $L a$ of $(Cayley(G_j,L),L \cdot 1)$, for all $a \in (G_j
\cap A) \setminus L$, with the $X_l$-monochromatic vertices $v_0
\cdot b$ of $C$, where $b \in G_l \cap A$ such that $b=_G a$. The
following lemma is a straightforward result of this construction.

\begin{lem} \label{glue_cayley graph to a vertex}
$(\Gamma',v_0')$ is a reduced precover of $G$ with
$Lab(\Gamma,v_0)=Lab(\Gamma',v_0')$, where $v_0$ is the basepoint
of $\Gamma$ and $v_0'$ is the (corresponding) basepoint of
$\Gamma'$.
\end{lem}
\begin{proof}
Obviously, by  construction, $\Gamma'$ is well-labelled,
compatible with $A$ and each monochromatic component of $\Gamma'$
is a cover of either $G_1$ or $G_2$. Thus $\Gamma'$ is a precover
of $G$. Moreover, $\Gamma'$ has no redundant components and
condition (2) from  Definition~\ref{def:reduced
 precover} is satisfied.   Hence $(\Gamma',v_0')$ is a reduced precover of $G$.

By  construction, $\Gamma$ and $Cayley(G_j,L,L \cdot 1)$ embed in
$\Gamma'$. Hence $(\Gamma,v_0) \subseteq (\Gamma',v_0')$, thus
$Loop(\Gamma,v_0) \subseteq Loop(\Gamma',v'_0)$. Therefore
$Lab(\Gamma, v_0) \subseteq Lab(\Gamma',v'_0)$.

 Let  $u \in Lab(\Gamma',v'_0)$. Hence there is $t'
\in Loop(\Gamma',v'_0)$ such that $lab(t')=_G u$.
If   $t'$ is a path in $\Gamma$ therefore $$lab(t')=_G u \in
Lab(\Gamma,v_0).$$
Otherwise there is a decomposition
$$t'=t'_1q_1t'_2q_2 \ldots q_{k-1}t'_k$$ such that
$\iota(t'_1)=\tau(t'_k)=v'_0$, and for all $1 \leq i \leq k$, $ \;
t_i' \subseteq \Gamma$ and $q_i$ is a path in $ \Gamma'$ which
doesn't exist in  $\Gamma$.

Thus for all $1 \leq i \leq k$, $q_i$ is a path in $
Cayley(G_j,L,L \cdot 1)$ such that $\iota(q_i)=v_{i_1}$ and
$\tau(q_i)=v_{i_2}$ are the common images in $\Gamma'$ of the
vertices $w_{i_1},w_{i_2}  \in \{v_0 \cdot a \; | \; a \in A
\setminus L \}$ of $C$ and the vertices $u_{i_1},u_{i_2} \in \{L a
\; | \; a \in A \setminus L \}$ of $ Cayley(G_j,L,L \cdot 1)$,
respectively.
By abuse of notation, we  write
$v_0 \cdot a_{i_1}=w_{i_1}=\iota(q_i) =u_{i_1}=L a_{i_1} \ {\rm
and} \ v_0 \cdot a_{i_2}=w_{i_2}=\tau(q_i) =u_{i_2}=L a_{i_2},$
where $a_{i_1}, \; a_{i_2} \in A \setminus L$.

Since $(L a_{i_1}) \cdot lab(q_i)=L a_{i_2}$, there exists $b \in
L$ such that $lab(q_i) =_G a_{i_1}^{-1} b a_{i_2}$. Hence $w_{i_1}
\cdot (a_{i_1}^{-1} b a_{i_2})= (v_0 \cdot b) \cdot a_{i_2}=v_0
\cdot a_{i_2}=w_{i_2},$ because $b \in L \leq K$.
Therefore  there exists a path $q_i'$  in $C$ (that is in
$\Gamma$) such that
$$\iota(q_i')=w_{i_1}, \ \tau(q_i')=w_{i_2}, \ lab(q_i')  =_G  lab(q_i) .$$
Thus there exists  a path $t$ in $\Gamma$ such that
 $ t=t'_1q'_1t'_2q'_2 \ldots q'_kt'_k$.  Therefore
\begin{eqnarray}
lab(t) & \equiv & lab(t'_1)lab(q'_1)lab(t'_2)lab(q'_2) \ldots
lab(q'_k)lab(t'_k) \nonumber \\
    &=_{G} & lab(t'_1)lab(q_1)lab(t'_2)lab(q_2) \ldots
lab(q_k)lab(t'_k) \nonumber \\
    &\equiv & lab(t'). \nonumber
\end{eqnarray}
 Since  $ lab(t') \in Lab(\Gamma, v_0)$, we have
$Lab(\Gamma)=Lab(\Gamma')$.

\end{proof}

%----------------------------------------------------------------------------------------

\begin{lem} \label{lem: normal path via vertices of red-precover}
Let $(\Gamma,v_0)$ be a reduced precover of $G$. Then for  each $v
\in V(\Gamma)$ there exists a normal path $p$ in $\Gamma$ closed
at $v_0$ such that $v \in V(p)$.
\end{lem}
\begin{proof}

Let $C$ be a $X_i$-monochromatic component of $\Gamma$ ($i \in
\{1,2\}$) such that $v \in V(C)$.

Since $C$ is not redundant, by Lemma \ref{red-precover=normal
paths}, there exists a normal path $q$ in $\Gamma$ closed at $v_0$
that goes through $C$. Let $q=q_1 \cdots q_m$ be a normal
decomposition of $q$ into maximal monochromatic paths. Assume that
$q_l \subseteq q \cap C$ ($1 \leq l \leq m$). Let $v_1=\iota(q_l)$
and $v_2=\tau(q_l)$.

If $v \in V(q_l)$ then  $p=q$ is the desired path.
Otherwise, we proceed in the following way. Assume first that
$m=1$. Then, by the proof of Lemma \ref{red-precover=normal
paths}, $Lab(C,v_0) \neq \{1\}$. Let $t$ be a path in $C$ with
$\iota(t)=v_0$, $\tau(t)=v$ and $lab(t) \equiv g$.  Hence
$Lab(C,v)=g^{-1}Lab(C, v_0)g \neq \{1\}$. Therefore there exists a
nonempty path $q' \in Loop(C,v)$ such that $lab(q') \neq_{G} 1$.
Therefore $tq'\overline{t} \in Loop(C,v_0)$ and
$lab(tq'\overline{t}) \neq_{G} 1$. Thus $lab(p_vp\overline{p_v})$
is a normal word, because it is a nonempty word of syllable length
1, which is not the identity in $G$. Hence $p=tq'\overline{t}$ is
the desired normal path in $\Gamma$ closed at $v_0$ that goes
through $v$.

Assume now that $m>1$. Let $t_j$ be paths in $C$ with $\iota(t_j)
= v_j$ to $\tau(t_j)=v$ ($j \in \{1,2\}$), see Figure \ref{figure:
NormalPathViaVertexOfPrecover}. Let $t=t_1\overline{t_2}$. Since
$deg_{\Gamma}(v) \geq 2$ ($\Gamma$ is a precover of $G$), we can
assume that $t$ is freely reduced.

%---------------------------------------------------

\begin{figure}[!htb]
\begin{center}
\psfrag{p1 }{$q_1$} \psfrag{p2 }{$q_{l-1}$}
\psfrag{pl }{$q_l$}
\psfrag{p3 }{$q_{l+1}$} \psfrag{pm }{$q_m$}
\psfrag{t1 }[][]{$t_1$} \psfrag{t2 }[][]{$t_2$}
\psfrag{C }[][]{$C$}
\psfrag{v0 }[][]{$v_0$} \psfrag{v }[][]{$v$}
\includegraphics[width=0.5\textwidth]{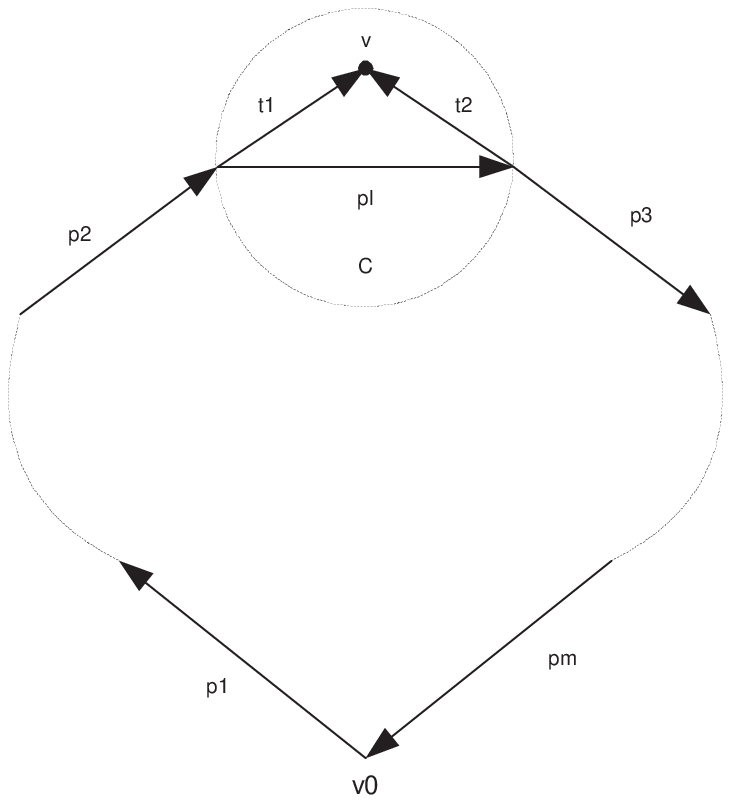}
\caption{ \label{figure: NormalPathViaVertexOfPrecover}}
\end{center}
\end{figure}

%---------------------------------------------------

If $lab(t) \not\in A$ then the path $p=q_1 \cdots q_{l-1}tq_{l+1}
\cdots q_m$ is the desired normal path in $\Gamma$ closed at $v_0$
which goes through $v$.

If $lab(t) \in A$ then $lab(t \overline{q_l}t ) \equiv lab(t)lab(
\overline{q_l})lab(t) \not\in A$, because $lab(q_l) \not\in A$.
Hence $p=q_1 \cdots q_{l-1}(t \overline{q_l}t )q_{l+1} \cdots q_m$
is the desired normal path in $\Gamma$ closed at $v_0$ which goes
through $v$.

\end{proof}

%%%%%%%%%%%%%%%%%%%%%%%%%%%%%%%%%%%%%%%%%%%%%%%%%%%%%%%%%%%%%%%%%%%%%%%%%%%%%%%%%%%%%%%%%%

\section{The Main Theorem} \label{sec: main theorem}

Let $H$ be a finitely generated subgroup of the amalgam  $G=G_1
\ast_A G_2$. As was discussed in the previous sections, there
exist labelled graphs which can be 'naturally' associated with
$H$. Thus the examples of such graphs are the normal core of
$Cayley(G,H)$, on the one hand, and a reduced precover of $G$,
$(\Gamma,v_0)$, with $Lab(\Gamma,v_0)=H$, on the other.

Below we prove that  normal cores and  reduced precovers
determining the same subgroup $H$ are the same. That is they
define the same part of $Cayley(G,H)$ in different ways: the
normal core defines it theoretically, while the reduced precover
characterizes it constructively.

%-------------------------------------------------------------------------------------------
\begin{thm}[The Main Theorem] \label{thm: unique reduced precover}
Let $H$ be a finitely generated subgroup of the amalgam  $G=G_1
\ast_A G_2$. Then up to isomorphism there exists a unique reduced
precover of $G$ determining $H$, which is the normal core
$(\Delta,H \cdot 1)$ of $Cayley(G,H)$.
\end{thm}

We separate the proof of the main theorem into two parts. First we
prove that if there exists  a reduced precover of $G$ determining
the subgroup $H$ then up to isomorphism it is unique. This
statement follows from  Theorem~\ref{isom_finite_grp}. Then we
prove (Theorem~\ref{ncore=reduced precover}) that given a finitely
generated subgroup $H$ of $G$ there exists  a reduced precover
determining $H$, which is precisely the normal core of
$Cayley(G,H)$.

%%%%%%%%%%%%%%%%%%%%%%%%%%%%%%%%%%%%%%%%%%%%%%%%%%%%%%%%%%%%%%%%%%%%
Let $(\Gamma,v_0)$ be a pointed graph labelled with $X^{\pm}$.
Define
$$\mu:(\Gamma,v_0) \rightarrow (Cayley(G,S),S \cdot 1)$$
 such that
$$\forall v \in V(\Gamma), \ \mu(v)=(S \cdot 1) \cdot lab(p)=S (lab(p)),$$
where $p$ is a path in $\Gamma$ with
 $\iota(p)=v_0$, $ \tau(p)=v $, and
$$\forall e \in
 E(\Gamma), \  \mu(e)=(\mu(\iota(e)),lab(e)).$$
In the proof of Lemma \ref{lemma1.5} (Lemma 1.5 in \cite{gi_sep})
Gitik shows that $\mu$ is a morphism of labelled pointed graphs
which is injective if $\Gamma$ is $G$-based. %
Hence if $\Gamma$ is a precover of $G$, then the morphism $\mu$ is
an embedding. We are interested in an identification of the
monomorphic image $\mu(\Gamma)$ inside $Cayley(G,S)$.

%%%%%%%%%%%%%%%%%%%%%%%%%%%%%%%%%%%%%%%%%%%%%%%%%%%%%%%%%%%%%%%%%%%%%%%%%%%%%%%%%%%%%%%%%%%%%%%%%%%%%%%%%%

%--------------------------------------------------------------------------------------------

\begin{thm} \label{isom_finite_grp}
Let $(\Gamma,v_0)$ be a  reduced precover of $G$.

Let $(\Delta, H \cdot 1)$ be the normal core of $Cayley(G,H)$,
where $H=Lab(\Gamma,v_0)$.

Then $\mu(\Gamma, H \cdot 1)=(\Delta, H \cdot 1)$.
\end{thm}
%------------------------------------------------------------------------
\begin{proof}
If $V(\Gamma)=\{v_0\}$ and $E(\Gamma)=\emptyset$ then $H=\{1\}$ by
Remark \ref{remark: trivial subgroup}. Therefore, by Remark
\ref{rm core=union of closed paths_finite grps}, $V(\Delta)=\{H
\cdot 1\}$ and $E(\Gamma)=\emptyset$. We are done.

%%%------------------------------------------------------------

First we show that $\mu(\Gamma,v_0) \subseteq (\Delta, H \cdot
1)$.
%--------------------------------------------------------------
%
Let $u \in V(\mu(\Gamma))=\mu(V(\Gamma))$. Hence $u=\mu(v)$, where
$v \in V(\Gamma)$. Without loss of generality, we can assume that
$v \neq v_0$, otherwise the statement is trivial ($\mu(v_0)=H
\cdot 1  \in V(\Delta)$), because $\mu$ is a morphism of pointed
graphs.

By Lemma \ref{lem: normal path via vertices of red-precover},
there exists a normal path $p$ in $\Gamma$ closed at $v_0$ such
that $v \in V(p)$.
Since graph morphisms commute with $\iota$, $\tau$ and preserve
labels, $\mu(p)$  is a normal path in $Cayley(G,H)$ that goes
through the vertex $\mu(v)=u$. Thus $\mu(p)$ is a  path in the
normal core $(\Delta,H \cdot 1)$ and $\mu(v)=u \in V(\Delta)$.
Therefore  $V(\mu(\Gamma)) \subseteq V(\Delta)$. Since graph
morphisms commute with $\iota$, $\tau$ and preserve labels, we
conclude  that $\mu(\Gamma,v_0) \subseteq (\Delta, H \cdot 1)$.

%------------------------------------------------------------------------------------------

Now we prove that $\mu(\Gamma,v_0) \supseteq (\Delta, H \cdot 1)$.
%----------------------------------------------------------

Let $\sigma \in V(\Delta)$. Then there is a normal path $\delta$
in $\Delta$ closed at $H \cdot 1$ in $Cayley(G,H)$ such that
$\sigma \in V(\delta)$. Thus $lab(\delta) \in H$ is a word in
normal form. Hence there exists a path $p$ in $\Gamma$ closed at
$v_0$  with $lab(p)=_G lab(\delta)$. Since $\Gamma$ is a precover,
by Lemma~\ref{lem: all normal paths are in precover}, there exists
a path $p'$ in $\Gamma$ closed at $v_0$  with $lab(p') \equiv
lab(\delta)$. Therefore $\delta=\mu(p')$. Hence there exists $v
\in V(p')$ such that $\sigma=\mu(v)$.

Therefore $V(\Delta) \subseteq V(\mu(\Gamma))$. Since graph
morphisms commute with $\iota$, $\tau$ and preserve labels, we
conclude that $(\Delta, H \cdot 1) \subseteq \mu(\Gamma,v_0)$.
Hence $(\Delta, H \cdot 1) = \mu(\Gamma,v_0)$.

\end{proof}

\begin{cor} \label{cor: mu=isomorphism}
Following the notation of Theorem \ref{isom_finite_grp}, $\mu$ is
an isomorphism of $(\Gamma, v_0)$ and $(\Delta, H \cdot 1)$.
\end{cor}

\begin{cor} \label{cor: red precovers of the same subgroup are
isomorphic}
Any pair of  reduced precovers of $G$ determining the same
subgroup are isomorphic.
\end{cor}
%\begin{proof}
%Let $(\Gamma_1,v_1)$ and $(\Gamma_2,v_2)$ be reduced precovers
%such that  $Lab(\Gamma_1,v_1)= Lab(\Gamma_2,v_2)=H \leq G$. Then,
%by Theorem \ref{isom_finite_grp}, $\mu_1(\Gamma_1,v_1) = (\Delta,H
%\cdot 1)=\mu_2(\Gamma_2,v_2)$, where $\mu_i: (\Gamma_i,v_i)
%\rightarrow (Cayley(G,H),SH \cdot 1)$ ($i \in \{1,2\}$) is the
%monomorphism defined above. By Corollary \ref{cor:
%mu=isomorphism}, $\mu_1\mu_2^{-1}$ is an isomorphism of
%$(\Gamma_1,v_1)$ and $(\Gamma_2,v_2)$.
%\end{proof}

%%%%%%%%%%%%%%%%%%%%%%%%%%%%%%%%%%%%%%%%%%%%%%%%%%%%%%%%%%%%%%%%%%%%%%%%%

\begin{thm} \label{ncore=reduced precover}
Let $H $ be a finitely generated subgroup of $G$. Then the normal
core $(\Delta,H \cdot 1)$ of $Cayley(G,H)$ is a reduced precover
of $G$ with $Lab(\Delta,H \cdot 1)=H$.
\end{thm}
\begin{proof}
Without loss of generality, we can assume that $H \neq \{1\}$,
because otherwise, by Remark \ref{rm core=union of closed
paths_finite grps}, the statement is trivial.

By  definition, a well-labelled graph $\Gamma$  is a precover of
$G$ if it is $G$-based and each $X_i$-monochromatic component of
$\Gamma$ ($i \in \{1,2\}$) is a cover of $G_i$.

Since $\Delta$ is a subgraph of $Cayley(G,H)$,  $\Delta$ is
well-labelled with $X_1^{\pm} \cup X_2^{\pm}$ and $G$-based.
Therefore each $X_i$-monochromatic component of $\Delta$ is
$G_i$-based ($i \in \{1,2\}$).  By Lemma~\ref{lemma1.5}, in order
to conclude that each such component is a cover of $G_i$, we have
to show that it is $X_i^{\pm}$-saturated.

Let $C$ be a $X_i$-monochromatic component of  $\Delta$ ($i \in
\{1,2\}$). Let $v \in V(C)$ and $x \in X_i$. Let $C'$ be the
$X_i$-monochromatic component of $Cayley(G,H)$ such that $C
\subseteq C'$. Therefore there is $e \in E(C')$ such that $lab(e)
\equiv x, \ \iota(e)=v$ and $v_x=\tau(e) \in V(C')$.

Since $v \in V(C) \subseteq V(\Delta)$,  there is a normal form
path $p$ in $\Delta$ that goes through $v$. If $e \in E(p)$ then
we are done. Otherwise, let
$$p=p_1 \cdots p_{l-1} q p_{l+1} \cdots p_k$$ be a normal
decomposition of $p$ into maximal monochromatic subpaths, such
that $p \cap C=q$, $v \in V(q)$ and $lab(q) \in G_i \setminus A$
(see Figure \ref{figure proof normal core=red precover}).
%-------------------------------------------------------------------------------------------

\begin{figure}[!htb]
\begin{center}
\psfrag{p1 }{$p_1$} \psfrag{p2 }{$p_{l-1}$}
\psfrag{p3 }{$p_{l+1}$} \psfrag{pk }{$p_k$}
\psfrag{q1 }[][]{$q_1$} \psfrag{q2 }[][]{$q_2$}
\psfrag{v }{$v$} \psfrag{vx }{$v_x$}
\psfrag{e }{$e$} \psfrag{t }{$t$} \psfrag{C }[][]{$C$}
\psfrag{v1 }[][]{$H \cdot 1$}
\includegraphics[width=0.5\textwidth]{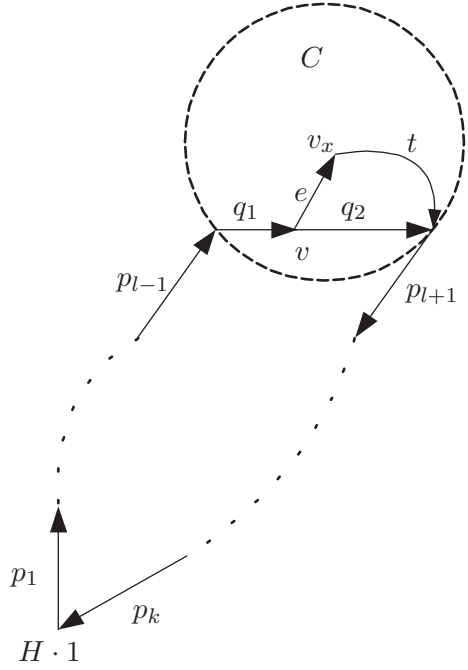}
\caption[An auxiliary figure for the proof of
Lemma~\ref{ncore}]{\footnotesize An auxiliary figure for the proof
of Lemma~\ref{ncore} \label{figure proof normal core=red
precover}}
\end{center}
\end{figure}

%-------------------------------------------------------------------------------------------

Let $\iota(q)=v_1, \ \tau(q)=v_2$, $q=q_1q_2$ such that $
\tau(q_1)=v=\iota(q_2)$.
Let $t$ be a path in $C'$  with $\iota(t)=v_x, \ \tau(t)=v_2$.
Thus $\iota(q_1et)=v_1=\iota(q)$ and $ \tau(q_1et)=v_2=\tau(q)$.
(Since the graph $Cayley(G,H)$ has no hairs, without loss of
generality, we can assume that the path $q_1et$ is freely
reduced.)

Let $$q'=\left\{%
\begin{array}{ll}
    q_1et, & \hbox{$lab(q_1et) \not\in A$;} \\
    (q_1et) \overline{q}(q_1et), & \hbox{$lab(q_1et) \in G_1 \cap A$.} \\
\end{array}%
\right.$$
Thus the path $q'$ has the same endpoints as $q$, $v_x \in V(q')$
and $lab(q') \in G_i \setminus A$. Hence $p=p_1 \cdots p_i q'
p_{i+1} \cdots p_k$ is a normal form path that goes through the
vertex $v_x$. Therefore $p \subseteq \Delta$ and $e \in E(C)
\subseteq E(\Delta)$. As this is true for every $x \in X_i$, the
vertex $v$ is $X_i^{\pm}$-saturated.
Hence, by Definition \ref{def: precover}, $\Delta$ is a precover
of $G$.

%----------------------------------------------------------------------------------------------

By Lemma \ref{red-precover=normal paths},  $\Delta$ has no
redundant monochromatic components, because for each $v \in
V(\Delta)$ there is a path in normal form closed at $H \cdot 1$
that goes through $v$.

Assume now that $C$ is a $X_i$-monochromatic component of $\Delta$
($i \in \{1,2\}$) such that $(C,H \cdot 1)$ is isomorphic to
$Cayley(G_i, K, K \cdot 1)$, where $K \cap A$ is a nontrivial
subgroup of A. Then there exists a nonempty normal path $p$ in $ C
\subseteq \Delta$ closed at $H \cdot 1$ with $lab(p) \equiv w \in
K \cap A \cap G_i$. Since $\{1\} \neq K \cap A \leq A$, there
exists $1 \neq u \in G_j \cap A$ ($1 \leq i \neq j \leq 2$) such
that $w=_G u$. Thus the syllable length of the words $w$ and $u$
is equal to $1$. Therefore these words are in normal form.

The graph $Cayley(G,H)$ is $X^{\pm}$ saturated and compatible with
$A$. Thus $H \cdot 1 \in VB(Cayley(G,H))$ and therefore there
exists a path $q$ in \linebreak[4] $Cayley(G,H)$ closed at $H
\cdot 1$ with $lab(q) \equiv u$. Hence $q \subseteq \Delta$,
because $u$ is in normal form. Since $\Delta$ is a precover of
$G$, $D \subseteq \Delta$, where $D$ is a $X_j$-monochromatic
component of $Cayley(G,H)$ such that $q \subseteq D$ and $(D,H
\cdot 1)$ is isomorphic to $Cayley(G_j, L, L \cdot 1)$.

Since,  $\Delta$ is compatible with $A$ (as a subgraph of
$Cayley(G,H)$), $L \cap A=_G K \cap A$. Then, by
Definition~\ref{def:reduced precover}, $(\Delta, H \cdot 1)$ is a
reduced precover of $G$.

\end{proof}

\begin{proof}[Proof of The Main Theorem]
The statement is an immediate consequence of Corollary \ref{cor:
red precovers of the same subgroup are isomorphic}, Theorem
\ref{ncore=reduced precover} and Lemma \ref{ncore}.

\end{proof}

%%%%%%%%%%%%%%%%%%%%%%%%%%%%%%%%%%%%%%%%%%%%%%%%%%%%%%%%%%%%%%%%%%%%%%%%%%%

\label{subsec:TheConstruction_finite_grp}

\section{The Algorithm} \label{subsec: TheAlgorithm}

Let $H$ be a finitely generated subgroup of an amalgam $G=G_1
\ast_A G_2$.   By Definition~\ref{def: normal core} and Remark
\ref{rm core=union of closed paths_finite grps}, the normal core
of $Cayley(G,H)$ depends on $H$ itself and not on the set of
subgroup generators, therefore this graph is canonically
associated with the subgroup $H$. Hence it can be exploited to
study certain  properties of $H$.

In  Lemma \ref{ncore} we prove that when the factors $G_1$ and
$G_2$ are finite groups,  the normal core of $Cayley(G,H)$ is a
finite graph, which is completely defined by $H$. Thus,
evidentially, it can be constructed. Our main theorem (Theorem
\ref{thm: unique reduced precover}) hints the way. Indeed, by
Theorem \ref{thm: unique reduced precover}, the normal core of
$Cayley(G,H)$ is the unique reduced precover of $G$ determining
$H$. Therefore in order to construct the normal core of
$Cayley(G,H)$ we should take the `right bunch' of copies of
relative Cayley graphs of the free factors, glue them to each
other according to the amalgamation, and verify that the obtained
precover is reduced. If not then it can be converted to a reduced
precover using Lemmas \ref{Cayley-remove} and \ref{glue_cayley
graph to a vertex}.

The precise algorithm, the proof of its finiteness and validity,
and the complexity analysis are presented in the current section.

%---------------------------------------------------------------------------------

Our proof of the finiteness of the normal core is based  on the
following result of Gitik \cite{gi_quas}.
\begin{defin}[\cite{gi_quas}] \label{def: projection map}
Let $G=gp\left\langle X|R \right\rangle$. Let $$\pi_S : Cayley(G)
\rightarrow Cayley(G,S)$$ be the  projection map such that
$\pi_S(g)=Sg$ and $\pi_S(g,x)=(Sg,x)$.

A \underline{geodesic} in $Cayley(G,S)$ is the image of a geodesic
in $Cayley(G)$ under the projection $\pi_S$. The
\underline{geodesic core} of $Cayley(G,S)$, \fbox{$Core(G,S)$}, is
the union of all closed geodesics in $Cayley(G,S)$ beginning at
the vertex $S \cdot 1$.
\end{defin}

\begin{lem}[Lemma 1.5 in \cite{gi_quas}] \label{lemma-quas}
A subgroup $S$ of a group $G$ is $K$-quasiconvex in $G$ if and
only  if $Core(G,S)$ belongs to the $K$-neighborhood of $S \cdot
1$ in $Cayley(G,S)$.
\end{lem}

\begin{lem} \label{ncore}
    Let $H$ be a finitely generated subgroup of $G=G_1 \ast_A G_2$.

    If $G_1$ and $G_2$ are finite groups.
    Then  the normal core $(\Delta, H \cdot 1)$ of $Cayley(G,H)$ is a finite
    graph.
\end{lem}
\begin{proof}
Since the group $G$ is locally-quasiconvex (\cite{ikap}), the
subgroup $H$ is quasiconvex. Therefore, $Core(G,H)$ is a finite
graph, by Lemma~\ref{lemma-quas}.

Let $\overline{\gamma}$ be a closed  normal path starting at $H
\cdot 1$ in $(\Delta,H \cdot 1) \subset (Cayley(G,H),H \cdot 1)$.
Thus $\overline{\gamma}$ is the image under the projection map
$\pi_H$ (see Definition~\ref{def: projection map}) of the normal
path $\gamma$ in $Cayley(G)$ whose endpoints and the label are in
$H$. That is $lab(\gamma) \equiv h \in H$.

Since $G_1$ and $G_2$ are finite, they are quasiconvex subgroups
of the hyperbolic group $G$. Thus  the conditions of
Lemma~\ref{lemma4.1} are satisfied. Let $\epsilon \geq 0$ be the
constant from Lemma~\ref{lemma4.1}. Let $\delta$ be  a geodesic in
$Cayley(G)$ with the same endpoints as $\gamma$. By
Lemma~\ref{lemma4.1},   there exists a strong normal path
$\delta'$ in $Cayley(G)$ with the same endpoints as $\delta$ such
that $\delta' \subset N_{\epsilon}(\delta)$ and $\delta \subset
N_{\epsilon}(\delta')$ (see Figure
\ref{fig:NormalPathsCloseToGeodesics}).

\begin{figure}[!htb]
\begin{center}
\psfragscanon
\psfrag{a }[][]{$< \epsilon$}
\includegraphics[width=0.7\textwidth]{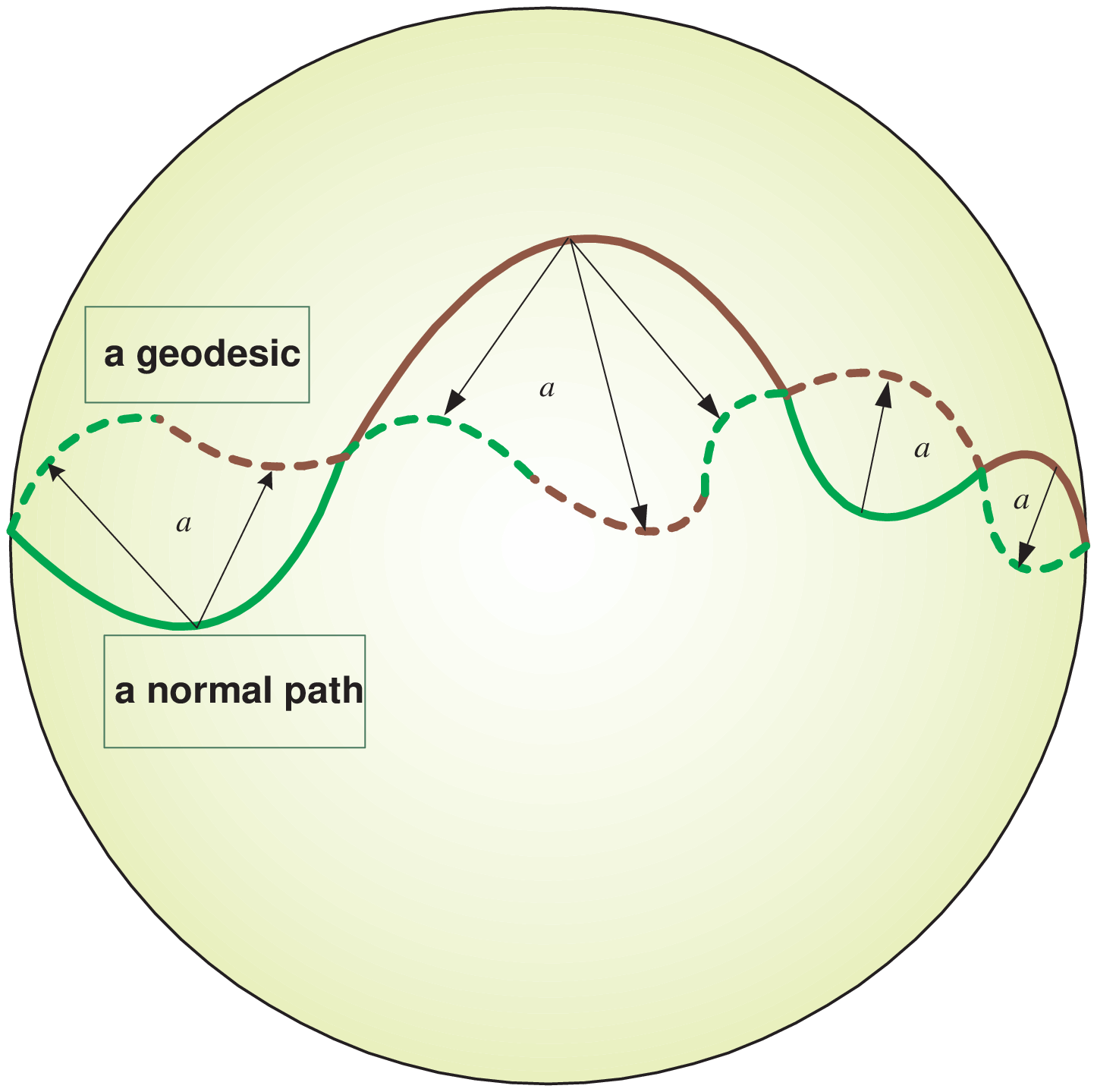}
\caption{ \label{fig:NormalPathsCloseToGeodesics}} \end{center}
\end{figure}

Thus $\gamma$ and $\delta'$ are two normal form  paths  in
$Cayley(G)$ with the same endpoints. Therefore $lab(\gamma)=_G
lab(\delta')$. By Corollary \ref{cor: distance between normal
paths}, $\gamma \subset N_{d}(\delta')$ and  $\delta' \subset
N_{d}(\gamma)$, where $d=max(diameter(G_1), diameter(G_2))$.

 Let $\epsilon'=\epsilon+d$. Then
$\gamma \subset N_{\epsilon'}(\delta)$ and  $\delta \subset
N_{\epsilon'}(\gamma)$. Since the projection map $\pi_H$ does not
increase distances,  and it maps $\gamma$ onto $\overline{\gamma}$
in $(\Delta,H \cdot 1)\subseteq (Cayley(G,H), H \cdot 1)$ and
$\delta$ onto $\overline{\delta}$ in $Core(G,H)\subseteq
(Cayley(G,H), H \cdot 1)$, we have $\overline{\gamma} \subset
N_{\epsilon'}(\overline{\delta})$ and $\overline{\delta} \subset
N_{\epsilon'}(\overline{\gamma})$.

This implies that $Core(G,H) \subset N_{\epsilon'}(\Delta)$ and
$\Delta \subset N_{\epsilon'}(Core(G,H))$. Since $Core(G,H)$ is a
finite graph we conclude that the graph $(\Delta, H \cdot 1)$ is
finite as well.

\end{proof}

Below we follow the notation of Grunschlag \cite{grunschlag},
distinguishing between the ``\emph{input}'' and the ``\emph{given
data}'', the information that can be used by the algorithm
\emph{``for free''}, that is it does not affect the complexity
issues.

\begin{center}
\large{\emph{\underline{\textbf{Algorithm}}}}
\end{center}

\begin{description}
\item[Given] Finite groups $G_1$, $G_2$, $A$ and the amalgam
$G=G_1 \ast_{A} G_2$ given via $(1.a)$, $(1.b)$ and $(1.c)$,
respectively.

We assume that the Cayley graphs and all the relative Cayley
graphs of the free factors are given.
\item[Input]  A finite set $\{ g_1, \cdots, g_n \} \subseteq G$.
\item[Output] A finite graph $\Gamma(H)$ with a basepoint $v_0$
which is a reduced precover of $G$ and the following holds
\begin{itemize}
 \item
$Lab(\Gamma(H),v_0)=_{G} H$;
 \item $H=\langle g_1, \cdots, g_n \rangle$;
 \item a normal word $w$ is in $H$ if and only if
  there is a loop (at $v_0$) in $\Gamma(H)$
labelled by the word $w$.
 \end{itemize}
\item[Notation] $\Gamma_i$ is the graph obtained after the
execution of the $i$-th step.

%
%------------------------------------------
%
\medskip

    \item[\underline{Step1}] Construct a based set of $n$ loops around a common distinguished
vertex $v_0$, each labelled by a generator of $H$;
%%%%
    \item[\underline{Step2}] Iteratively fold edges and cut hairs;
%%%%
%   \item[\underline{Step3}] Iteratively identify endpoints of each path labelled with a relator.
%    After each identification, if necessary, fold edges and "cut hairs".
%%%%%
 \item[\underline{Step3}] { \ }\\
%----
\texttt{For} { \ } each $X_i$-monochromatic component $C$ of
$\Gamma_2$ ($i=1,2$) { \ } \texttt{Do} \\
\texttt{Begin}\\
    pick an edge $e \in E(C)$; \\
    glue a copy  of $Cayley(G_i)$   on $e$ via identifying $ 1_{G_i} $  with $\iota(e)$ \\
    and identifying the two copies of $e$ in $Cayley(G_i)$ and in $\Gamma_2$; \\
    \texttt{If}  { \ } necessary  { \ } \texttt{Then} { \ } iteratively fold
    edges; \\
 \texttt{End;}

%%%%%
 \item[\underline{Step4}] { \ } \\
%----
\texttt{ For}  { \ } each $v \in VB(\Gamma_3)$ { \ } \texttt{ Do} \\
 \texttt{If} { \ } there are paths $p_1$ and $p_2$, with $\iota(p_1)=\iota(p_2)=v$
 and $\tau(p_1)~\neq~\tau(p_2)$  such that
 $$lab(p_i) \in G_i \cap A \ (i=1,2) \ {\rm and} \  lab(p_1)=_G
 lab(p_2)$$
\texttt{ Then} { \ } identify $\tau(p_1)$ with $\tau(p_2)$; \\
 \texttt{If}  { \ } necessary  { \ } \texttt{Then} { \ } iteratively fold
    edges; \\
%%%%%

 \item[\underline{Step5}]
%----
%
Reduce  $\Gamma_4$ by iteratively removing all \emph{redundant}
 $X_i$-monochromatic components $C$ such that
\begin{itemize}
 \item $(C,\vartheta)$ is isomorphic to $Cayley(G_i, K, K \cdot 1)$, where $K \leq A$ and
$\vartheta \in VB(C)$;
 \item  $|VB(C)|=[A:K]$;
 \item one of the following holds
    \begin{itemize}
        \item  $K=\{1\}$ and $v_0 \not\in VM_i(C)$;
        \item $K$ is  a nontrivial subgroup of $A$ and $v_0  \not\in V(C)$.\\
    \end{itemize}
\end{itemize}

Let $\Gamma$ be the resulting graph;\\

\texttt{If}  { \ }
$VB(\Gamma)=\emptyset$ and $(\Gamma,v_0)$ is isomorphic to $Cayley(G_i, 1_{G_i})$ \\
\texttt{Then} { \ } we set $V(\Gamma_5)=\{v_0\}$ and
$E(\Gamma_5)=\emptyset$; \\
\texttt{Else} { \ } we set $\Gamma_5=\Gamma$.

%%%%%
 \item[\underline{Step6}] { \ } \\
%----
 \texttt{If} { \ }
 \begin{itemize}
  \item $v_0 \in VM_i(\Gamma_5)$ ($i \in \{1,2\}$);
  \item $(C,v_0)$ is isomorphic to $Cayley(G_i,K,K \cdot 1)$, where $L=K \cap A$ is a nontrivial
  subgroup of
 $A$ and $C$ is a $X_i$-monochromatic component of $\Gamma_5$ such that $v_0 \in V(C)$;
  \end{itemize}
\texttt{Then} { \ } glue to $\Gamma_5$ a $X_j$-monochromatic
component ($1 \leq i \neq j \leq 2$) $D=Cayley(G_j,L,L \cdot 1)$
via identifying $L \cdot 1$ with $v_0$ and \\
identifying the vertices $L \cdot a$ of  $Cayley(G_j,L,L \cdot 1)$
with the vertices $v_0 \cdot a$ of $C$, for all $a \in A \setminus
L$.

%%%%
Denote $\Gamma(H)=\Gamma_6$.

\end{description}
%-------------------------------------------------------------------------------------------

%-------------------------------------------------------------

\begin{remark} \label{stal-mar-meak-kap-m}
{\rm The first two steps of the above algorithm correspond
precisely to the Stallings' folding algorithm for finitely
generated subgroups of free groups (see \cite{stal, mar_meak,
kap-m}). This allows one to  refer to our algorithm as the
\emph{generalized Stallings' (folding) algorithm} for finitely
generated subgroups of amalgams of finite groups.

By the results of \cite{stal, mar_meak, kap-m}, the graph
$\Gamma_2$ is finite, well-labelled with $X^{\pm}$, has no hairs
and $Lab_{F(X)}(\Gamma_2,v_0)=H$, where $Lab_{F(X)}(\Gamma_2,v_0)$
is the image of $lab(Loop(\Gamma_2,v_0))$ in the free group
$F(X)$.

} \e
\end{remark}

%------------------------------------------------------------

\begin{figure}[!h]
\psfrag{x }[][]{$x$} \psfrag{y }[][]{$y$} \psfrag{v }[][]{$v$}
\psfrag{x1 - monochromatic vertex }[][]{{\footnotesize
$\{x\}$-monochromatic vertex}}
\psfrag{y1 - monochromatic vertex }[][]{\footnotesize
{$\{y\}$-monochromatic vertex}}
\psfrag{ bichromatic vertex }[][]{\footnotesize {bichromatic
vertex}}
\includegraphics[width=\textwidth]{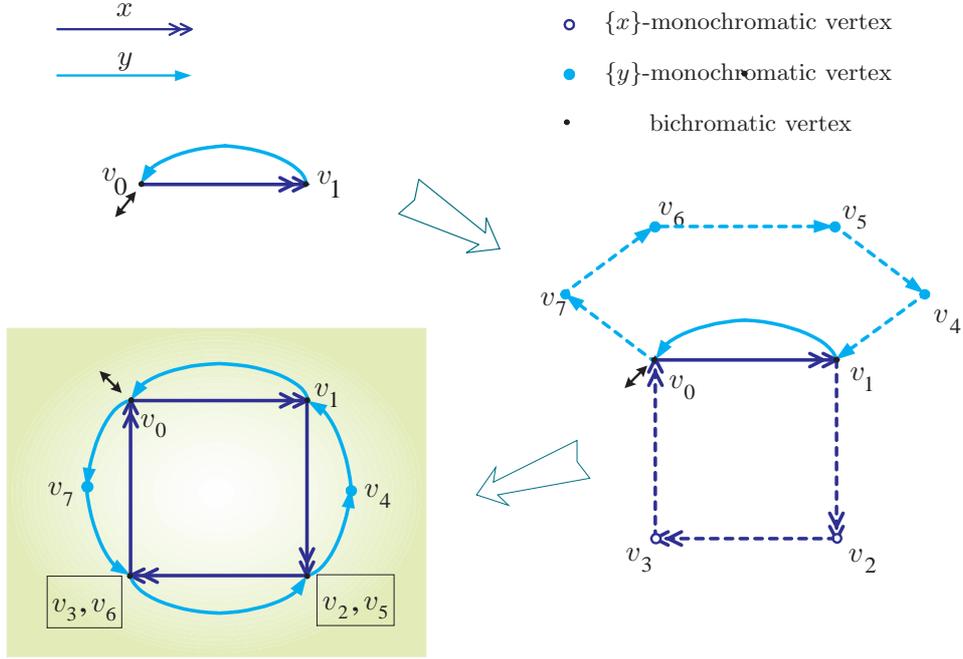}
\caption[The construction of $\Gamma(H_1)$]{ \footnotesize {The
construction of $\Gamma(H_1)$.}
 \label{fig: example of H=xy}}
\end{figure}

\begin{figure}[!htb]
\psfrag{v }[][]{$v$}
\includegraphics[width=\textwidth]{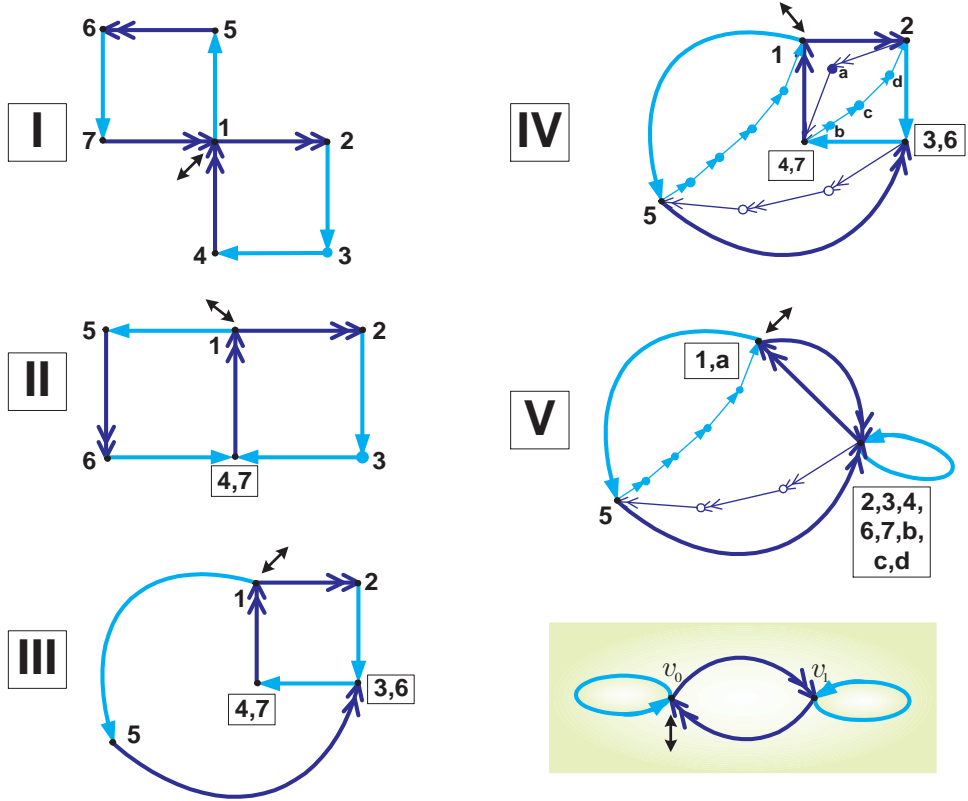}
\caption[The construction of $\Gamma(H_2)$]{ \footnotesize {The
construction of $\Gamma(H_2)$.} \label{fig: example of H=xy^2x,
yxyx}}
\end{figure}

\begin{ex} \label{example: graphconstruction}
{\rm Let $G=gp\langle x,y | x^4, y^6, x^2=y^3 \rangle$.

Let $H_1$ and $H_2$ be  finitely generated subgroups of $G$ such
that
$$H_1=\langle xy \rangle \ {\rm and} \ H_2=\langle xy^2, yxyx \rangle.$$

The construction of $\Gamma(H_1)$ and $\Gamma(H_2)$ by the
algorithm presented above is illustrated on Figure \ref{fig:
example of H=xy}
 and  Figure  \ref{fig: example of H=xy^2x, yxyx}.}
\e
\end{ex}

\begin{lem} \label{construction_finite_grp}
  The  algorithm terminates and constructs
  the graph $(\Gamma(H),v_0)$ which is a finite reduced precover of $G$ with  $Lab(\Gamma(H),v_0)=H$.
\end{lem}
\begin{proof}
By Remark \ref{stal-mar-meak-kap-m}, the first two steps of the
algorithm terminates and construct the finite graph $\Gamma_2$.
Since $G_1$ and $G_2$ are finite groups, $Cayley(G_1)$ and
$Cayley(G_2)$ are finite graphs. Therefore, by the construction,
all the intermediate graphs $\Gamma_i$ ($ 3 \leq i \leq 6$) are
finite. Moreover they are constructed by a finite sequence of
iterations. Thus the resulting graph $\Gamma(H)$ is finite.

By Remark \ref{stal-mar-meak-kap-m} and by
Lemma~\ref{foldings-cutting-hairs},
$Lab(\Gamma_2,v_0)=Lab_{F(X)}(\Gamma_2,v_0)=H$.

Applying to each of the intermediate graphs $\Gamma_i$ ($ 3 \leq i
\leq 6$) the appropriate lemma from
Lemmas~\ref{foldings-cutting-hairs}, \ref{glue_cayley graph to an
edge} (see Appendix), \ref{identificacion-def-relator},
\ref{Cayley-remove} and \ref{glue_cayley graph to a vertex}, we
get
 $$Lab(\Gamma_6,v_0)=Lab(\Gamma_5,v_0)=Lab(\Gamma_4,v_0)=Lab(\Gamma_3,v_0)=Lab(\Gamma_2,v_0)=H.$$
 Thus $Lab(\Gamma(H),v_0)=H$.

Graphs $\Gamma_3$ and $\Gamma_4$ are well-labelled with $X^{\pm}$,
due to the folding operations, by \cite{stal}.

$\Gamma_3$ has no hairs. Indeed, since the graphs $\Gamma_2$ and
$Cayley(G_i)$  ($i \in \{1,2\}$) have no hairs, the intermediate
graph of the third step obtained after the gluing operations has
no hairs. Moreover, the graphs $\Gamma_2$ and $Cayley(G_i)$  ($i
\in \{1,2\}$) are well-labelled. Thus the only possible foldings
in the intermediate graph are between edges of a
$X_i$-monochromatic component $C$ of $\Gamma_2$ and edges of the
copy of $Cayley(G_i)$  ($i \in \{1,2\}$) glued to $C$ along the
common edge $e$. Therefore  the terminal vertices of the resulting
edges  have degree greater than 1.

Since foldings keep $X_i^{\pm}$-saturated vertices
$X_i^{\pm}$-saturated and keep closed paths closed, the image of a
copy of $Cayley(G_i)$ ($i\in\{1,2\}$) in $\Gamma_3$ remains
$G_i$-based and $X_i^{\pm}$-saturated. Thus by
Lemma~\ref{lemma1.5}, it is a cover of $G_i$.

 Let $C$ be a
$X_i$-monochromatic component of $\Gamma_2$ from the definition of
the third step, that is $e \in E(C)$. Let $C'$ be its image in
$\Gamma_3$. Then $C' \subseteq S$, where $S$ is an image of a copy
of the Cayley graph $Cayley(G_i)$ ($i\in\{1,2\}$) in $\Gamma_3$.
Indeed, let $v' \in V(C')$ be the image of the vertex $ v \in
V(C)$. Hence, since $C$ is connected, there exist a path $q$ in
$C$ such that $\iota(q)=\iota(e)$ and $\tau(q)=v$. Thus $lab(q)
\in (X_i^{\pm})^*$. Since the graph operations of the third step
can be viewed as graph morphisms, they preserves labels and
``commutes'' with endpoints. Thus the image $  q'$ of the path $q$
in $C'$ satisfies $\iota(q')=\vartheta$, $\tau(q')=v'$ and
$lab(q') \equiv lab(q)$, where $\vartheta$ is the ``common'' image
in $\Gamma_3$ of the vertices $\iota(e)$ of $\Gamma_2$ and
 $1_{G_i}$ of $Cayley(G_i)$.

On the other hand, since $Cayley(G_i)$ is $X_i^{\pm}$ saturated,
there exists a path $\gamma$ in $Cayley(G_i)$ with
$\iota(\gamma)=1_{G_i}$ and $lab(\gamma) \equiv lab(q)$. Then
there is a path $\gamma'$ in $S$ with $\iota(\gamma')=\vartheta$
and $lab(\gamma') \equiv lab(\gamma) \equiv lab(q)$. Since
$\Gamma_3$ is well-labelled, we have $q'=\gamma'$. Hence $V(C')
\subseteq V(S)$. Thus $C' \subseteq S$.

Therefore all $X_i$-monochromatic components of $\Gamma_3$  are
covers of $G_i$ ($i\in\{1,2\}$).

Let $v \in VB(\Gamma_3)$ and let $p_1$ and $p_2$ be paths in
$\Gamma_3$ such that
\begin{itemize}
 \item $\iota(p_1)=\iota(p_2)=v$;
 \item  $\tau(p_1) \neq \tau(p_2)$;
 \item $lab(p_i) \in G_i \cap A \ (i=1,2)$;
 \item $ lab(p_1)=_G lab(p_2)$.
\end{itemize}
Then $v \in VB(\Gamma_3)$ and $deg(v) \geq 2$.

Let $\nu$ be a vertex  which is the result of the identification
of the vertices $\tau(p_1)$ and $\tau(p_2)$ of $\Gamma_3$. If
$\tau(p_1)$ and $\tau(p_2)$ are monochromatic
 vertices of $\Gamma_3$ of different colors, then no
foldings are possible at $\nu$ and $deg(\nu) \geq 2$.

Otherwise at least one of them is bichromatic in $\Gamma_3$. Then
$\nu$ is a bichromatic vertex of $\Gamma_4$ and  foldings are
possible at $\nu$. However foldings keep bichromatic vertices
bichromatic. Thus $\nu \in VB(\Gamma_4)$ and $deg(\nu) \geq 2$.

Therefore, since $\Gamma_3$ has no hairs,  $\Gamma_4$ has no hairs
as well. By Remarks \ref{remark: identif-relator} and \ref{remark:
morphism of precovers}, each $X_i$-monochromatic component of
$\Gamma_4$ is a cover of $G_i$ ($i \in \{1,2\}$).
By the construction, $\Gamma_4$  is compatible.  Hence, by
Lemma~\ref{corol2.13}, the graph $\Gamma_4$ is a precover of $G$.

By Lemma~\ref{Cayley-remove}, $\Gamma_5$ is a precover of $G$ as
well. Since, by the construction, $\Gamma_5$   has no redundant
monochromatic components, Lemma \ref{glue_cayley graph to a
vertex} implies that $\Gamma(H)=\Gamma_6$ is a reduced precover of
$G$.

\end{proof}

Now we sketch the complexity analysis of the above algorithm.

\begin{lem}[Complexity] \label{complexity of construction}
Let $m$  be the sum of the lengths of words $g_1, \ldots g_n$.
Then the algorithm computes $(\Gamma(H),v_0)$ in time $O(m^2)$.
\end{lem}
\begin{proof} As is well known, see \cite{b-m-m-w}, the construction of the
bouquet can be done in time proportional to $m$, foldings can be
implemented in time proportional to $m^2$ and cutting hairs can be
done in time proportional to $m$. Therefore the first two steps of
the algorithm are completed in time $O(m^2)$, while the graph
$\Gamma_2$ satisfies: $|E(\Gamma_2)| \leq m$ and $|V(\Gamma_2)|
\leq m$.

Given Cayley graphs of both free factors $G_1$ and $G_2$, the
gluing operations of the third step take time proportional to $m$,
because we just identify one edge of each monochromatic component
of $\Gamma_2$ (whose detecting takes $|E(\Gamma_2)|$) with the
corresponding edge of the graph $Cayley(G_i)$, $i \in \{1,2\}$.

 Note that
$|V(\Gamma_3)|=k_1 \cdot |G_1| +k_2 \cdot |G_2|$, where $k_i$, $i
\in \{1,2\}$, is a number of $X_i$-monochromatic components of
$\Gamma_2$. Since  the information about the factors $G_1$ and
$G_2$ is given, that is it is not a part of the input, and since
$k_1+k_2 \leq m$, we conclude that the number $|V(\Gamma_3)|$ is
proportional to $m$. Similarly, $|E(\Gamma_3)|$ is proportional to
$m$ as well.

The detecting of bichromatic vertices of $\Gamma_3$ takes time
proportional to $|V(\Gamma_3)|$, that is it takes time
proportional to $m$. By the proof of
Lemma~\ref{construction_finite_grp} (the proof of the fourth
step), there are at most $|A|$ identifications for each
bichromatic vertex of $\Gamma_3$. Thus the identifications of the
fourth step take at most  $|VB(\Gamma_3)| \cdot |A|$. However, the
description of the third step implies that
$$|VB(\Gamma_3)|=|VB(\Gamma_2)| \leq
|V(\Gamma_2)|.$$
Since the number of vertices of the intermediate graph of the
fourth step obtained just after the above identifications is not
greater  than $|V(\Gamma_3)|$, the foldings operations applied to
this graph can be implemented in time proportional to
$|V(\Gamma_3)|^2$, by \cite{b-m-m-w}. Since $|V(\Gamma_3)|$ is
proportional to $m$, it takes time proportional to $m^2$. Thus,
summarizing the analysis of the fourth step, we see that its
implementation takes
 $O(m^2)$.

The indication of connected monochromatic components of $\Gamma_4$
takes time proportional to $|E(\Gamma_4)|$. Since $|E(\Gamma_4)|
\leq |E(\Gamma_3)|$ and $|E(\Gamma_3)|$ is proportional to $m$,
this procedure takes time proportional to $m$.

By the proof of Lemma~\ref{construction_finite_grp}, the graph
$\Gamma_4$ is a precover of $G$, hence its $X_i$-monochromatic
components are covers of $G_i$ for all $i \in \{1,2\}$.  Since the
information about the factors $G_1$ and $G_2$ is given, that is it
is not a part of the input, the verifications concerning
monochromatic components of $\Gamma_4$ take $O(1)$.

%Let $C$ be a monochromatic component of $\Gamma_4$. Detecting its
%bichromatic vertices takes $|V(C)| \leq |G_i|$ for $i \in
%\{1,2\}$, while verifying whether it is isomorphic to the relative
%Cayley graph of a subgroup $K$ of $A$ takes time ($|K|!$).

Since in the worst case the monochromatic component of $\Gamma_4$
that has to be deleted via the fifth step might appear at the end
of the verification process, while it induces a series of
deletions, the fifth step can be completed  in time proportional
to $|E(\Gamma_4)|$, that is in $O(m^2)$.

The last step of the algorithm takes at most $|A|$, that is
constant according to our assumption (it is a part of the `given
information'').

Summarizing the above description of the steps complexity, we
conclude that the algorithm constructs the resulting graph
$\Gamma(H)$ in time  $O(m^2)$.

\end{proof}

\begin{remark} \label{rem:construction_exponential}
{\rm  Note that if the group presentations of the free factors
$G_1$ and $G_2$, as well as the monomorphisms  between the
amalgamated subgroup $A$ and the free factors are a part of the
input (the \emph{uniform version} of the algorithm) then we have
to build the groups $G_1$ and $G_2$ (that is to construct their
Cayley graphs and relative Cayley graphs).

Since we assume that the groups $G_1$ and $G_2$ are finite,  the
Todd-Coxeter algorithm and the Knuth Bendix algorithm are suitable
\cite{l_s, sims, stil} for these purposes. Then the complexity of
the construction depends on the group presentation of $G_1$ and
$G_2$ we have: it could be even exponential in the size of the
presentation. Therefore the generalized Stallings' folding
algorithm with these additional constructions could take time
exponential in the size of the input.}

\e
\end{remark}

%Moreover, by Lemma~\ref{complexity of construction}, this
%algorithm is polynomial, when the information about the free
%factors is given, that is it is not a part of the input. This
%condition can be viewed as a disadvantage of our algorithm,
%because its removing turns the algorithm to be exponential.

%%%%%%%%%%%%%%%%%%%%%%%%%%%%%%%%%%%%%%%%%%%%%%%%%%%%%%%%%%%%%%%%%%%%%%%%%%%%%%%%%%%%%%%%%%%%%%%%%%%%%%%%%%%

%\section{Normal Core}
%\label{sec:Normal Core}

\begin{thm} \label{thm amalgam of finite grps}
Let $Y$ be a finite subset of $G$ and let $H=\langle Y \rangle$ be
a finitely generated subgroup of $G$. Then the resulting graph
$(\Gamma(H),v_0)$ constructed by the generalized Stallings'
folding algorithm  is the normal core of $Cayley(G,H)$.
\end{thm}
\begin{proof} The generalized Stallings' folding algorithm constructs a graph $(\Gamma(H),v_0)$, which is a
finite reduced precover of $G$ with $Lab(\Gamma(H),v_0)=H$, by
Lemma~\ref{construction_finite_grp}. Hence, by
Theorem~\ref{isom_finite_grp}, $(\Gamma(H),v_0)$ is isomorphic to
$(\Delta,H \cdot 1)$, the normal core of $Cayley(G,H)$. Since this
isomorphism is unique, by Remark \ref{unique isomorphism}, the
graph $(\Gamma(H),v_0)$ can be identified with the normal core of
$Cayley(G,H)$.
\end{proof}

\begin{remark}[Canonicity and Constructibility] \label{rem: ncore is constructible}
{\rm   Theorem~\ref{thm amalgam of finite grps} implies that the
normal core of a relative Cayley graph is constructible.

Since, by Definition~\ref{def: normal core} and Remark \ref{rm
core=union of closed paths_finite grps}, the normal core of
$Cayley(G,H)$ depends on $H$ itself and not on the set of subgroup
generators, Theorem~\ref{thm amalgam of finite grps} implies that
the graph $(\Gamma(H),v_0)$ is canonically associated with $H$.
 } \e
\end{remark}

As an immediate consequence of Theorem~\ref{thm amalgam of finite
grps} we get the following corollary, which provide a solution for
the membership problem for finitely generated subgroups of
amalgams of finite groups. We discuss it in the next section.

\begin{cor}\label{remark: normal elements in H}
A  normal word $g$ is in $H$ if and only if it labels a closed
path in $\Gamma(H)$ starting at $v_0$.
\end{cor}
\begin{proof}
A normal word $g$ is in $H$ if and only if it labels a normal path
in the normal core of $Cayley(G,H)$ closed at $H \cdot 1$. Since,
by Theorem~\ref{thm amalgam of finite grps}, $(\Gamma(H),v_0)$
constructed by the generalized Stallings' folding algorithm  is
the normal core of $Cayley(G,H)$, we obtain the desired
conclusion.

\end{proof}

%%%%%%%%%%%%%%%%%%%%%%%%%%%%%%%%%%%%%%%%%%%%%%%%%%%%%%%%%%%%%%%%%%%%%%%%%%

\section{The Membership Problem} \label{sec: membership pr}

The \emph{membership problem} (or the \emph{generalized word
problem}) for a subgroup of a given group asks to decide whether a
word in the generators of the group is an element of the given
subgroup.

As is well known (\cite{a_g}), the membership problem for finitely
generated subgroups  is solvable in amalgams of finite groups.
Different types of solutions can be found in \cite{c-otto,
holt-hurt, grunschlag} and other sources.

Below we introduce a solution of the membership problem for
finitely generated subgroups of amalgams of finite groups which
employs subgroup graphs (normal cores) constructed by the
generalized Stallings' foldings algorithm, presented in
Section~\ref{subsec: TheAlgorithm}.

\begin{cor} \label{membership_finite_grp}
Let $g, h_1, \ldots h_n \in G$.
Then there exists an algorithm which  decides whether or not $g$
belongs to the subgroup $H=\langle h_1, \ldots, h_n \rangle$~of~
$G$.
\end{cor}
\begin{proof} First we construct the graph $\Gamma(H)$, using the
algorithm from Section~\ref{subsec: TheAlgorithm}. By
Corollary~\ref{remark: normal elements in H},  $g \in H$ if and
only if there is a normal path $p$ in $\Gamma(H)$ closed at the
basepoint $v_0$ such that $lab(p) =_G g$. That is the word
$lab(p)$ is a normal form of the word $g$.

Thus in order to decide if $g \in H$ we have to begin with a
calculation of a normal form $\overline{g}$ of the given word $g$.
If $g$ is a normal word then we just skip the calculation and put
$\overline{g} \equiv g$ .
Otherwise we use a well-known rewriting procedure \cite{l_s} to
find $\overline{g}$. This usage is possible because the membership
problem for the amalgamated subgroup $A$ is solvable in the free
factors $G_1$ and $G_2$ (indeed, they  are finite groups).

Now we have to verify if there exists a path $p$ in $\Gamma(H)$
closed at the basepoint $v_0$ such that $lab(p) \equiv
\overline{g}$.
It can be done as follows. We start at the vertex $v_0$ and try to
read the word $\overline{g}$ in the graph $\Gamma(H)$. If we
become stuck during this process or if we don't return to the
vertex $v_0$ at the end of the word $\overline{g}$, then $g$ is
not in $H$. Otherwise we conclude that $g \in H$.

\end{proof}

\begin{ex}
{\rm  Let $H_2$ be the subgroup of $G$ from Example \ref{example:
graphconstruction}. Then using Figure \ref{fig: example of
H=xy^2x, yxyx} and the algorithm described in Corollary
\ref{membership_finite_grp}, we easily conclude that $xyx \in
H_2$, because $v_0 \cdot (xyx)=v_0$ in $\Gamma(H_2)$. But
$xy^3x^{-5} \not\in H_2$, because $v_0 \cdot (xy^3x^{-5}) \neq
v_0$. } \e
\end{ex}

The algorithm presented with the proof of Corollary
\ref{membership_finite_grp} provides a solution for the
\emph{membership problem} for finitely generated subgroups of
amalgams of finite groups with the following description.
\begin{description}
    \item[GIVEN] Finite groups $G_1$, $G_2$, $A$ and the amalgam
$G=G_1 \ast_{A} G_2$ given via $(1.a)$, $(1.b)$ and $(1.c)$,
respectively.

We assume that the Cayley graphs and all the relative Cayley
graphs of the free factors are given.
    \item[INPUT]  Words $ g, h_1, \ldots,
h_n \; \in \; G$.
    \item[DECIDE]  Whether or not $g$ belongs to the subgroup $H=\langle h_1, \ldots,
h_n \rangle$.
\end{description}

\subsubsection*{Complexity.} {Let $m$ be the sum of the lengths of
the words $h_1, \ldots h_n$. By Lemma~\ref{complexity of
construction}, the algorithm from Section \ref{subsec:
TheAlgorithm} computes $(\Gamma(H),v_0)$ in time $O(m^2)$. The
verification of the normality of the word $g$ is proportional to
$|g|$ and the computation of its normal form takes time
$O(|g|^{2})$. To read a normal word in the graph $(\Gamma(H),v_0)$
in the way, explained in the proof of Corollary
\ref{membership_finite_grp}, takes time equal to the length of the
word.
Therefore the complexity of the algorithm is $O(m^2+|g|^2)$.}

If in the above description the input is changed to:
\begin{description}
\item[INPUT]  Words $ h_1, \ldots, h_n \; \in \; G$ and a normal
word $g \in G$.
\end{description}
then the complexity of the algorithm will be $O(m^2+|g|)$.

In some papers, the following slightly different description of
the \emph{membership problem} can be found.
\begin{description}
    \item[GIVEN]  Finite groups $G_1$, $G_2$, $A$ and the amalgam
$G=G_1 \ast_{A} G_2$ given via $(1.a)$, $(1.b)$ and $(1.c)$,
respectively.

We assume that the Cayley graphs and all the relative Cayley
graphs of the free factors are given.

     The subgroup $H=\langle h_1, \ldots,
h_n \rangle$ of $G$.
    \item[INPUT]  A normal word $g \in G$.
    \item[DECIDE]  Whether or not $g$ belongs to the subgroup $H$.
\end{description}

 In this context  the subgroup $H$ is given, that is $(\Gamma(H),v_0)$
is constructed and can be used for free. Therefore the complexity
of this  algorithm is linear in the length of the word $g$,
because we simply have  to read it in the graph $(\Gamma(H),v_0)$
which takes time equal to $|g|$.

%Comparing this analysis to the result of Holt and Hurt
%(\cite{holt-hurt}), one can see that our algorithm for the
%solution of the membership problem  given by the last description
%is faster. Recall that given a coset automatic system, the
%algorithm presented in \cite{holt-hurt} solves the generalized
%word problem for the subgroup $H$ of a group $G$ in time $O(n^2)$,
%where $n$ is the length of the input word, while our algorithm
%solves it in time $O(n)$. Note that if the input word $g$ is not
%normal then our algorithms have the same the complexity.

Another variation of the \emph{membership problem} is the
\emph{uniform membership problem}, when the presentation of the
group $G$ is a part of the input.

\begin{description}
    \item[GIVEN]   -
    \item[INPUT]    Finite groups $G_1$, $G_2$, $A$ and the amalgam
$G=G_1 \ast_{A} G_2$ given via $(1.a)$, $(1.b)$ and $(1.c)$,
respectively.

Words $ g, h_1, \ldots, h_n \; \in \; G$.
    \item[DECIDE]  Whether or not $g$ belongs to the subgroup $H=\langle h_1, \ldots,
h_n \rangle$.
\end{description}

\subsubsection*{Complexity.}{ The algorithm given along with the
proof of Corollary \ref{membership_finite_grp} provide a solution
for the above problem. However now the complexity of the
construction of $(\Gamma(H),v_0)$ might be exponential in the size
of the presentation $(1.a)$-$(1.c)$, by Remark
\ref{rem:construction_exponential}. Therefore the complexity of
the algorithm might be exponential in the size of the input.}

%%%%%%%%%%%%%%%%%%%%%%%%%%%%%%%%%%%%%%%%%%%%%%%%%%%%%%%%%%%%%%%%%%%%%%%%%%%%%%%%%%

\appendix
\label{appendix}

\section{}

When constructing graphs for subgroups of non free groups,
nontrivial relations of these groups have to be taken into
account. Roughly speaking, they have to be ``sewed'' somehow on
subgroup graphs. This gives rise to  a \emph{gluing} operation on
the graph. Moreover, if one is interested to construct a precover
of an amalgam  then the gluing operation of copies of Cayley
graphs of the free factors to the graph have to be defined.

Let $\Gamma$ be a graph well-labelled with $X^{\pm}$.  Let $e \in
E(\Gamma)$ such that $lab(e) \equiv w \in X_j$ ($j \in \{1,2\}$).

Let $\Gamma'$ be the graph constructed by taking the disjoint
union of the graphs $\Gamma$ and $Cayley(G_j)$ via the
identification of   the edge $e \in E(\Gamma)$ with  the edge $f
\in E(Cayley(G_j))$ such that $\iota(f)=1_{G_j}$ and $lab(f)
\equiv lab(e) $.

More precisely,
$$V(\Gamma')=(V(\Gamma) \setminus \{\iota(e),\tau(e)\}) \cup (V(Cayley(G_j)) \setminus \{\iota(f),\tau(f)\})
 \cup \{v_1, v_2\}.$$
$$E(\Gamma')=(E(\Gamma) \setminus \{ e \}) \cup (E(Cayley(G_j)) \setminus \{f\})
 \cup \{ l \}.$$

The endpoints and arrows for the edges of $\Gamma'$ are defined in
a natural way.
$$\iota_{\Gamma'}(\xi)=\left\{%
\begin{array}{ll}
    \iota_{\Gamma}(\xi), & \hbox{if $\xi \in E(\Gamma) \setminus \{ e \}$;} \\
    \iota_{Cayley(G_j)}(\xi), & \hbox{if $\xi \in E(Cayley(G_j)) \setminus \{f\}$;} \\
    v_1, & \hbox{if $\xi=l$;} \\
    v_2, & \hbox{if $\xi=\overline{l}$.} \\
\end{array}%
\right.    $$
We define labels on the edges of $\Gamma'$ as follows:
$$lab_{\Gamma'}(\xi) \equiv \left\{%
\begin{array}{ll}
    lab_{\Gamma}(\xi), & \hbox{if $\xi \in E(\Gamma) \setminus \{ e \}$;} \\
    lab_{Cayley(G_j)}(\xi), & \hbox{if $\xi \in E(Cayley(G_j)) \setminus \{f\}$;} \\
    w, & \hbox{if $\xi=l$.} \\
\end{array}%
\right.    $$

We say that $\Gamma'$ is obtained from $\Gamma$ by \emph{gluing a
copy of $Cayley(G_j)$ along the edge $e$ of $\Gamma$}.

\begin{lem} \label{glue_cayley graph to an edge}
Let $\Gamma'$ be the graph obtained from the well-labelled graph
$\Gamma$ gluing a copy of $Cayley(G_j)$ along the edge $e$ of
$\Gamma$.

Then $Lab(\Gamma,v_0)=Lab(\Gamma',v_0')$, where $v_0$ is the
basepoint of $\Gamma$ and $v_0'$ is the (corresponding) basepoint
of $\Gamma'$.
\end{lem}
\begin{proof}
Since the graph $Cayley(G_j)$  is $X_j^{\pm}$-saturated, there
exists an edge $f \in E(Cayley(G_j))$ such that $\iota(f)=1_{G_j}$
and $lab(f) \equiv lab(e)$. Thus the construction of $\Gamma'$ is
possible.

By the construction, $\Gamma$ and $Cayley(G_j)$ embed in
$\Gamma'$. Hence $(\Gamma,v_0) \subseteq (\Gamma',v_0')$, thus
$Loop(\Gamma,v_0) \subseteq Loop(\Gamma',v'_0)$. Therefore
$Lab(\Gamma, v_0) \subseteq Lab(\Gamma',v'_0)$.

 Let  $u \in Lab(\Gamma',v'_0)$. Then there is $t'
\in Loop(\Gamma',v'_0)$ such that \linebreak[4] $lab(t')=_G u$.
If  $t' \subseteq \Gamma$  then $lab(t')=_G u \in
Lab(\Gamma,v_0).$
Otherwise there is a decomposition
$$t'=t'_1q_1t'_2q_2 \ldots q_{k-1}t'_k,$$ such that
$\iota(t'_1)=\tau(t'_k)=v'_0$, and for all $1 \leq i \leq k$, $ \;
t_i'\subseteq \Gamma$ and $q_i$ is a path in $ \Gamma'$ which
doesn't exist in $\Gamma$.

Thus for all $1 \leq i \leq k$, $q_i$ is in $ Cayley(G_j)$ and
$\iota(q_i), \tau(q_i) \in \{v_1,v_2\}$, where $v_1$ and $v_2$ are
images of the vertices $\iota(e)$ and $\tau(e)$ of $\Gamma$ in
$\Gamma'$, respectively. Then either $\iota(q_i) = \tau(q_i)$ or
$\iota(q_i) \neq \tau(q_i)$. In the first case $q_i$ is a closed
path in $Cayley(G_j)$, hence $lab(q_i)=_{G_j} 1$. In the second
case either $\iota(q_i)=\iota(f)$, $\tau(q_i)=\tau(f)$ or
$\iota(q_i)=\tau(f)$, $\tau(q_i)=\iota(f)$. Therefore $lab(q_i)
=_{G_j} lab(f)$ or $lab(q_i) =_{G_j} (lab(f))^{-1}$, respectively.

 Let $ t=t'_1q'_1t'_2q'_2 \ldots q'_kt'_k$ be a  path in $\Gamma$ such
 that for all $1 \leq i \leq k$,
$$q_i'=\left\{%
\begin{array}{ll}
    \emptyset, & \hbox{$\iota(q_i) = \tau(q_i)$;}\\
    e, & \hbox{$lab(q_i)=_{G_j} lab(f)$;} \\
    \overline{e}, & \hbox{$lab(q_i) =_{G_j} (lab(f))^{-1}$.} \\
\end{array}%
\right. $$

By $q_i'=\emptyset$ we mean that $q_i'$ is the empty path, that is
$lab(q_i')=_{G_j}1$, with the desired initial-terminal vertex
$\iota(g_i')=\tau(q_i')=\tau(t_i')=\iota(t_{i+1}')$.

Since $lab(e) \equiv lab(f)$,
$$lab(q_i')=\left\{%
\begin{array}{ll}
    1, & \hbox{$\iota(q_i) = \tau(q_i)$;}\\
    lab(f), & \hbox{$lab(q_i)=_{G_j} lab(f)$;} \\
    (lab(f))^{-1}, & \hbox{$lab(q_i) =_{G_j} (lab(f))^{-1}$.} \\
\end{array}%
\right. $$ Thus $lab(q'_i)=_{G_j} lab(q_i)$. Therefore
\begin{eqnarray}
lab(t) & \equiv & lab(t'_1)lab(q'_1)lab(t'_2)lab(q'_2) \ldots
lab(q'_k)lab(t'_k) \nonumber \\
    &=_{G_j} & lab(t'_1)lab(q_1)lab(t'_2)lab(q_2) \ldots
lab(q_k)lab(t'_k) \nonumber \\
    &\equiv & lab(t'). \nonumber
\end{eqnarray}
 Since  $ lab(t') \in Lab(\Gamma, v_0)$, we have
$Lab(\Gamma)=Lab(\Gamma')$.

\end{proof}

%%%%%%%%%%%%%%%%%%%%%%%%%%%%%%%%%%%%%%%%%%%%%%%%%%%%%%%%%%%%%%%%%%%%%%%%%%%%%%%%%%%

\end{document}